\documentclass[12pt]{article}
\makeindex

\usepackage[latin1]{inputenc} \usepackage[T1]{fontenc} \usepackage{lmodern}
\usepackage{pgf}
\usepackage{epsfig,graphicx,color,amsmath,a4wide,amssymb,amsfonts,amsbsy,xy,latexsym,epsf,pstricks,verbatim,times}

\usepackage{pgf,color}
\usepackage{pgfarrows}

\definecolor{LightGrey}{rgb}{.85,.85,.85}
\definecolor{DarkGrey}{rgb}{.5,.5,.5}

\definecolor{Blue}{rgb}{.0,.0,0.9}
\definecolor{LightBlue1}{rgb}{.2,.4,0.9}
\definecolor{LightBlue2}{rgb}{.3,.5,0.9}
\definecolor{LightBlue3}{rgb}{.4,.6,0.9}
\definecolor{LightBlue4}{rgb}{.5,.7,.9}
\definecolor{LightBlue5}{rgb}{.6,.8,.9}
\definecolor{LightBlue6}{rgb}{.7,.9,.9}

\definecolor{Red}{rgb}{.9,.0,.0}
\definecolor{LightRed1}{rgb}{0.9,.2,.4}
\definecolor{LightRed2}{rgb}{0.9,.3,.5}
\definecolor{LightRed3}{rgb}{0.9,.4,.6}
\definecolor{LightRed4}{rgb}{.9,.5,.7}
\definecolor{LightRed5}{rgb}{.9,.6,.8}
\definecolor{LightRed6}{rgb}{.9,.7,.9}

\def\va{{\vec{\alpha}}}
\def\vg{{\vec{\gamma}}}
\def\vxr{{\vec{x}_R}}
\def\vur{{\vec{u}_R}}
\def\vcur{{\overrightarrow{u_R}^{(-1)}}}
\def\viota{{\vec{\iota}}}

\def\vb{{\vec{\beta}}}

\def\Car{\mathop{\rm Char }}
\def\conju{\mathop{\rm conj }}
\definecolor{Grey}{rgb}{.5,.5,.5}

\definecolor{Blue}{rgb}{.0,.0,0.9}
\definecolor{LightBlue1}{rgb}{.2,.4,0.9}
\definecolor{LightBlue2}{rgb}{.3,.5,0.9}
\definecolor{LightBlue3}{rgb}{.4,.6,0.9}
\definecolor{LightBlue4}{rgb}{.5,.7,.9}
\definecolor{LightBlue5}{rgb}{.6,.8,.9}
\definecolor{LightBlue6}{rgb}{.7,.9,.9}

\definecolor{Red}{rgb}{.9,.0,.0}
\definecolor{LightRed1}{rgb}{0.9,.2,.4}
\definecolor{LightRed2}{rgb}{0.9,.3,.5}
\definecolor{LightRed3}{rgb}{0.9,.4,.6}
\definecolor{LightRed4}{rgb}{.9,.5,.7}
\definecolor{LightRed5}{rgb}{.9,.6,.8}
\definecolor{LightRed6}{rgb}{.9,.7,.9}

\xyoption{all}
\usepackage[english]{babel}

\newcounter{noalgo}[section]
\newdimen\indentalgo
\newdimen\indentalgodec\indentalgo=0.0mm\indentalgodec=10mm
\renewcommand{\thenoalgo}{\thesection.\arabic{noalgo}}
\newcommand{\Ret}{{\bf return }}
\newcommand{\If}{\advance\indentalgo by \indentalgodec {\bf if }}
\newcommand{\Then}{{\bf then }}
\newcommand{\For}{\global\advance\indentalgo by \indentalgodec {\bf for }}
\newcommand{\To}{{\bf to }}
\newcommand{\By}{{\bf by }}
\newcommand{\Do}{{\bf do }}
\newcommand{\Endindent}{\global\advance\indentalgo by -\indentalgodec}

\newdimen\decalage \decalage=0.5cm
\newcounter{algo} 
\newcommand{\nna}{\vskip 0.05cm \hskip 0.25\decalage}
\newcommand{\nns}{\vskip 0.05cm \hskip 0.25\decalage \stepcounter{algo} {\bf \thealgo.}\hskip\indentalgo}

\let\set\mathbb
\def\<<{\leavevmode
  \raise0.28ex\hbox{$\scriptscriptstyle\langle\!\langle$}\nobreak
  \hskip -.6pt plus.3pt minus.2pt\,}
\def\>>{\,\nobreak\hskip -.6pt plus.3pt minus.2pt
  \raise0.28ex\hbox{$\scriptscriptstyle\rangle\!\rangle$}}

\def\Tr{\mathop{\rm{Tr}}\nolimits }

\def\CC{{\set C}}

\def\End{\mathop{\rm End }}
\def\FF{{\set F}}

\def\FQ{{\FF _Q}}

\def\Fq{{\FF _q}}
\def\Fqb{{\bar{\FF}_q}}

\def\bK{{\bf K  }}
\def\bL{{\bf L  }}

\def\QQ{{\set Q}}

\def\RR{{\set R}}

\def\UU{{\set U}}
\def\ZZ{{\set Z}}
\def\agot{{\mathfrak a}}
\def\tgot{{\mathfrak t}}

\def\bgot{{\mathfrak b}}

\def\cL{{\cal L}}

\def\cO{{\cal O}}

\def\cR{{\cal R}}
\def\cS{{\cal S}}

\def\cgot{{\mathfrak  c}}

\def\lgot{{\mathfrak l}}

\def\mmu{{\set \mu}}

\newtheorem{lemma}{Lemma}

\newtheorem{theorem}{Theorem}
\newtheorem{definition}{Definition}

\newenvironment{myproof}[1][\myproofname]{\par
  \normalfont \topsep6pt\relax
  \trivlist
\item[\hskip\labelsep
  \itshape
  #1.]\ignorespaces
}{%
  \endtrivlist\hfill$\square$
}
\providecommand{\myproofname}{Proof}

\let\origItemize=\itemize
\def\NoSpacing{
  \itemsep=0pt
  \parskip=0pt
  \parsep=0pt
  \partopsep=0pt
  \topsep=0pt
}
\renewenvironment{itemize}{\origItemize\NoSpacing}{\endlist}

\begin{document}
\author{Jean-Marc Couveignes\thanks{Institut de  Math\'ematiques  de Toulouse,
Universit\'e de Toulouse et CNRS, D{\'e}partement de Math{\'e}matiques et
Informatique, Universit{\'e} Toulouse 2, 5 all{\'e}es Antonio Machado, 31058 Toulouse
c{\'e}dex 9.} and Reynald
Lercier\thanks{\textsc{DGA}/\textsc{C\'ELAR}, La Roche Marguerite, F-35174 Bruz.}~\thanks{IRMAR, Universit\'e de Rennes 1, Campus de Beaulieu, F-35042 Rennes.}}\title{Elliptic periods for finite
fields\thanks{Research supported by the
    French D{\'e}l{\'e}gation G{\'e}n{\'e}rale pour l'Armement,
Centre d'\'Electronique de l'Armement and by the Agence Nationale de la
Recherche (projet blanc ALGOL).}}

\maketitle

\bibliographystyle{plain}

\begin{abstract}
We construct two  new families
of basis for finite field extensions.
Bases in the first
family, the so-called \textit{elliptic bases},  are not quite normal  bases, but  they allow
very fast Frobenius exponentiation while preserving sparse
multiplication formulas.
Bases in the second family, the so-called \textit{normal elliptic bases} are
normal bases and allow fast 
(quasi-linear) arithmetic. 
We prove  that all extensions admit models of this kind.
\end{abstract}


\section{Introduction}

The main computational advantage of normal basis for a
finite field extension $\FF_{q^d}/\Fq$ is that they allow fast
exponentiation by $q$ since it corresponds
to a cyclic shift of coordinates, and it can be computed in time
$O(d)$.
There is a  concern however
about how difficult is multiplication in this
context. 

Let $\alpha$ and $\beta$ be two elements in $\FF_{q^d}$ with
coordinates $\va = (\alpha_i)_{0\leqslant i\leqslant d-1}$ and $\vb =
(\beta_i)_{0\leqslant i \leqslant d-1}$ in the given normal basis. Let
$(\gamma_i)_{0\leqslant i\leqslant d-1}$ be the coordinates of the product $\alpha
\times \beta$. Each $\gamma_i$ is a bilinear form in $\va$ and $\vb$.
The number of non-zero terms in $\gamma_i$ does not depend on $i$
because the $d$ corresponding tensors are cyclic shifts of each
others.  This number of terms is called the {\it complexity} ${\mathcal C}$ of
the normal basis. 
Multiplication
with the straightforward  algorithm can be done
with $2d{\mathcal C}$
operations ($d{\mathcal C}$
when coefficients of the bilinear forms $\gamma_i$ are all $\pm 1$).
It was shown by Mullin, Onyszchuk, Vanstone and
Wilson  \cite{MOVW} that the complexity ${\mathcal C}$
is at least $2d-1$.
This bound is reached by the so-called optimal normal bases.  But such optimal normal bases only exist
for very special extensions.  As a general fact, normal bases with
bounded complexity are not known to exist, unless the degree $d$ takes
very special and sparse values.

Normal bases with low complexity usually are constructed using {\it
  Gauss periods} as in work by Ash, Blake and Vanstone \cite{AVS}
or Gao and Lenstra \cite{GL}.
The construction uses $r$-th roots
of unity where $r=kd+1$ is prime.  It requires that $q$ generates the
unique quotient of order $d$ of $(\ZZ/r\ZZ)^*$. The parameter $k$ is
very important and should be kept as small as possible, because the
complexity of the normal basis is bounded by $(d-1)k+d$  and is not expected
to be much smaller \cite[Theorem
4.1.4]{Gao}.  Optimal normal
bases occur when $k=1$ or $k=2$.  This corresponds to very sparse
values of $d$. In general, for $q$ a prime, assuming the Extended
Riemann Hypothesis, it has been shown by Adleman and Lenstra 
\cite{AL} that there exists a $k$ and a $r$ as above with
$r=O(d^4(\log (dq))^2)$. This is unfortunately of no use when bounding
the complexity.  In some cases, there is no $k$ at all  \cite[Satz
3.3.4]{Was}.  We shall not survey all the variants  and improvements
for this  method. We  just quote works
by Christopoulou, Garefalakis, Panario and Thomson
 \cite{CGPT} where traces 
of optimal normal bases
are shown to have a reasonable complexity in some  
special cases.   Wan and Zhou show  \cite{WZ} that
the dual
of  type I  optimal normal bases have
good complexity too.

Gao, von zur Gathen and Panario  show  \cite{GGP} that  fast 
multiplication
methods (like FFT) can be adapted to  normal bases constructed
with Gauss periods. They give a multiplication algorithm
in such a normal basis with complexity $O(dk\log(dk)\log|\log(dk)|)$.
This is a considerable progress for Gauss normal bases
with bounded $k$. But in the general case, $k$ being 
only upperbounded
by $O(d^3(\log (dq))^2)$, this     is just too large.

In his thesis \cite{Gao} Gao presented a new way of constructing 
normal bases with low complexity.
In Gao's construction, the Lucas torus and
its isogenies play an important, though implicit, role. 
Gao thus constructs  more normal bases with low complexity.
In our work, we consider the remaining
algebraic groups  of dimension one: elliptic curves. 
Since there are many elliptic curves, we can enlarge significantly
the number  of cases where  a normal basis with
fast multiplication exists.
\medskip

In order to state our results, we shall need the following definition
where $v_\ell$ stands for the valuation associated to  the prime
$\ell$.
\begin{definition}\label{definition:dq}
Let $p$ be a prime and $q$ a power of $p$. Let $d\geqslant 2$ be an integer.

We denote by  $d_q$  the unique positive integer such that for
every prime $\ell$
\begin{itemize}
\item $v_\ell(d_q)=v_\ell(d)$ if $\ell$ is prime to $q-1$,
\item $v_\ell(d_q)=0$ if $v_\ell(d)=0$,
\item $v_\ell(d_q)=\max(2v_\ell(q-1)+1, 2v_\ell(d))$ if $\ell$
divides both $q-1$ and $d$.
\end{itemize}
\end{definition}
For example, if $d=14$ and $q=654323$
then $q-1=2.19.67.257$ and $d_q=2^{3}.7$.

Note that $d_q=d$ whenever $d$ is prime to $q-1$.
\medskip

We now can state our first result.
\begin{theorem}\label{th:slow}

To  every 
couple $(q,d)$ with  $q$
a prime power and $d\geqslant 2$ an integer and
$d_q\leqslant q^{\frac{1}{2}},$
 one can associate  a normal basis $\Theta(q,d)$ of 
the degree $d$ extension of $\Fq$ such that the following
holds:
\begin{itemize}
\item There exist a positive constant $K$
and an algorithm that multiplies two elements given in 
the basis $\Theta(q,d)$ at  the expense of $5d^2+2d$ multiplications and $5d^2+4d$
  additions/subtractions in $\Fq$. The amount of necessary memory is
  $\leqslant Kd\log q$ bits.
\end{itemize}
\end{theorem}

There is also a fast arithmetic version of Theorem~\ref{th:slow}.
\begin{theorem}\label{th:fast}
To every 
couple $(q,d)$ with $q$ a prime power and
$d\ge 2$ an integer 
and $d_q \leqslant q^{\frac{1}{2}}$, one can associate
 a normal basis $\Theta(q,d)$ of 
the degree $d$ extension of $\Fq$ such that the following holds:
\begin{itemize}
\item There exist a positive constant 
$K$ and an algorithm that multiplies two elements
given in   the  basis $\Theta(q,d)$ at the expense of 
$Kd\log d\log |\log d|$
operations
in $\Fq$. 
\item There exists  an algorithm that divides two elements
given in the  basis $\Theta(q,d)$ at the expense of 
$$Kd(\log d)^2\log |\log d|$$
operations
in $\Fq$. 
\end{itemize}
\end{theorem}

The basis $\Theta(q,d)$ that appears in Theorem~\ref{th:slow} and
Theorem~\ref{th:fast} has a multiplication tensor that 
mainly consists of  $5$ convolution products.
We also construct a basis $\Omega(q,d)$ having a sparse
multiplication tensor. Sparsity is useful when using
such constrained  devices as circuits. Further, this basis
$\Omega(q,d)$ allows a faster elementary multiplication algorithm
than $\Theta(q,d)$. It is not quite a normal basis but
exponentiation by $q$ is still done  in 
linear time.
\begin{theorem}\label{th:omega}
To  every 
couple $(q,d)$ with $q$ a prime power and $d\ge 2$
an integer and
 $d_q \leqslant 2 q^{\frac{1}{2}}$, one can associate  a  basis $\Omega(q,d)$ of 
the degree $d$ extension of $\Fq$ such that the following holds: 
\begin{itemize}
\item 
There exist   a positive constant $K$ and
an algorithm that computes the $q$-th power 
of an element given 
in  basis $\Omega(q,d)$ at the expense of $d-1$  multiplications
and $2d-3$ additions in $\Fq$. The
  amount of necessary memory is $ \leqslant Kd\log q$ bits.
\item There exists
 an algorithm that multiplies two elements given in basis
  $\Omega(q,d)$  at the expense of $(31d^2+6d)/12$ multiplications, $d^2/12$
  inverses and $(37d^2+30d)/12$ additions/subtrac\-tions in $\Fq$. The
  amount of necessary memory is $ \leqslant Kd\log q$ bits.
\end{itemize}
\end{theorem}

The following  result  is valid 
 without any restriction.
\begin{theorem}\label{th:strong}
To  every 
couple $(q,d)$,   one can associate   a  model  $\Xi(q,d)$ of 
the degree $d$ extension of $\Fq$ such that the following holds :

There exists a  positive constant $K$ 
such that the following is true :

\begin{itemize}
\item Elements in $\FF_{q^d}$
are represented by vectors with less than $K d(\log  d)^2(\log(\log d))^2$
components in  $\Fq$.
\item Addition (resp. substraction)
 of two elements in $\FF_{q^d}$
requires less than $$K d(\log  d)^2(\log(\log d))^2$$ additions (resp. substractions)
in $\Fq$.
\item Exponentiation by $q$ consists in a circular
shift of the the   coordinates.
\item There exists an algorithm
that multiplies two elements at the expense of 
 $$K d(\log  d)^3|\log(\log d)|^3$$
multiplications/additions/substractions in $\Fq$.
\item There exists an algorithm
that divides two elements at the expense of 
 $$K d(\log  d)^4|\log(\log d)|^3$$
multiplications/additions/substractions in $\Fq$.
\end{itemize}
\end{theorem}

So, for every finite field extension, there exists a model that
allows both  fast
multiplication and fast application of the  Frobenius automorphism.
\medskip

In Section~\ref{section:relations}, we recall simple relations between
low degree elliptic functions.  We show in
Section~\ref{section:ellipticbasis} that evaluation of such functions
at a well chosen divisor produces an almost normal basis for the
residue field. Relations between elliptic functions result in nice
multiplication formulas in this basis.  Such bases have similar
properties to those constructed by Gao in his thesis: they have low
complexity.  This is shown in Subsection~\ref{para:complex}.  In
Section~\ref{section:fast}, we construct normal bases allowing fast
(quasi-linear) multiplication.  We show in
Section~\ref{section:beyond} that an elliptic basis exists for any degree
$d$ extension of $\Fq$ provided $d$ is not too large.  
We explain in Subsection~\ref{para:trans} what to
do when $d$ is large.  
In Subsection~\ref{para:fastinv},  we introduce
a polynomial basis that can be related 
efficiently to the elliptic (normal) basis.
We deduce a fast inversion algorithm
for elliptic normal bases.

We further support our claims with extensive experiments using the
computational algebra system \textsc{magma}~\cite{magma08}. We developed
for this task a package, named \textsc{ellbasis}, the sources of which are
available on the web page of the second author.

 \medskip

\noindent
{\bf Acknowledgments}: 
We thank C{\'e}cile Dartyge, Guillaume Hanrot,  Gerald Tenenbaum and Jie Wu
for pointing Iwaniec's result on Jacobsthal's problem  to us.

\section{Linear and quadratic
relations among elliptic functions}\label{section:relations}

In this section, we study the simplest elliptic functions: those
with degree $2$. We prove simple linear and quadratic relations 
between these functions.
The monography \cite{Silverman} by J. Silverman contains 
all the necessary background about elliptic curves. 
\medskip

Let  $\bK$ be a field and let $E$ be an elliptic curve over $\bK$. 
We assume $E$ is given by some Weierstrass equation
\begin{displaymath}
Y^2Z+a_1XYZ+a_3YZ^2=X^3+a_2X^2Z+a_4XZ^2+a_6Z^3\,.
\end{displaymath}
We set $x=X/Z$,
$y=Y/Z$ and $z=-x/y=-X/Y$, and we find 
\begin{eqnarray*}
x&=&\frac{1}{z^2}-\frac{a_1}{z}-a_2-a_3z+O(z^2)\,,\\
y&=&-\frac{1}{z^3}+\frac{a_1}{z^2}+\frac{a_2}{z}+a_3+O(z)\,.
\end{eqnarray*}
The involution $P=(x,y)\mapsto -P=(x,-y-a_1x-a_3)$ transforms
$z$ into 
$$z(-P)=\frac{x}{y+a_1x+a_3}=-z-a_1z^2-a_1^2z^3-(a_1^3+a_3)z^4+O(z^5)\,.$$

If $A$ is a geometric point on $E$, we denote by 
$\tau_A$ the translation by $A$. 
We denote by $z_A = z \circ \tau_{-A}$ the composition
of $z$ with the translation by $-A$. 
We define  $x_A$ and $y_A$ in a similar way.
The composition of $z_A$ with the involution fixing $A$
is
$-z_A-a_1z_A^2-a_1^2z_A^3-(a_1^3+a_3)z_A^4+O(z_A^5).$
The composition of $1/z_A$ with the involution fixing $A$
is
$-1/z_A+a_1 +a_3z_A^2+O(z_A^3).$
\medskip

If  $A$ and $B$ are two distinct geometric points on $E$, we denote
by $u_{A,B}$ the function  on $E$ defined as 
\begin{equation*}\label{eq:slope}
u_{A,B}=\frac{y_A-y(A-B)}{x_A-x(A-B)}\,.
\end{equation*}
It has  polar divisor $-[A]-[B]$.
 It is invariant by the
involution exchanging $A$ and $B$,
$$u_{A,B}(A+B-P)=u_{A,B}(P)\,.$$ 
Its Taylor expansion at $A$ is
\begin{math}
u_{A,B}=-{1}/{z_A}-x_{A}(B)z_A+(y_{A}(B)+a_3){z_A}^2+O(z_A^3)\,.
\end{math}
\medskip

If $C$ is any third geometric point, we set
$\Gamma(A,B,C)=u_{A,B}(C)$. This is the slope of the secant (resp. tangent) to
$E$  going through 
$C-A$ and $A-B$. It is well defined for any three points $A$, $B$, $C$ such
that $\# \{A,B,C\}\geqslant 2$. It is finite if and only if 
$\# \{A,B,C\} = 3$. 
We check 
\begin{equation}\label{eq:moins}
\Gamma(-A, -B, -C)=-\Gamma(A,B,C)-a_1.
\end{equation}

The  Taylor expansions of $u_{A,B}$ at $A$ and $B$ are
\begin{eqnarray*}
u_{A,B}&=&-\frac{1}{z_A}-x_{A}(B)z_A+(y_{A}(B)+a_3)z_A^2+O(z_A^3)\\
&=&
\frac{1}{z_B}-a_1+x_{A}(B)z_B+(y_{A}(B)+a_1x_A(B))z_B^2+O(z_B^3).
\end{eqnarray*}
As a consequence $u_{B,A}=-u_{A,B}-a_1$, $x_B(A)=x_A(B)$ and
$y_B(A)=-y_A(B)-a_1x_A(B)-a_3$ and
examination of Taylor expansions at $A$, $B$ and $C$ shows that
\begin{equation}\label{eq:somme}
u_{A,B}+u_{B,C}+u_{C,A}=\Gamma(A,B,C) -a_1
\end{equation}
and 
\begin{equation}\label{eq:sym}
\Gamma(A,B,C)=u_{B,C}(A)=u_{C,A}(B)=u_{A,B}(C)=-u_{B,A}(C)-a_1.
\end{equation}
We deduce
\begin{equation*}
u_{B,C}=u_{B,C}(A)-(x_{A}(C)-x_{A}(B))z_A+(y_{A}(C)-y_{A}(B))z_A^2+ O(z_A^3).
\end{equation*}
\medskip

By comparison of Taylor expansions at $A$, $B$ and $C$
we prove
\begin{displaymath}
u_{A,B}u_{A,C}=x_A+u_{B,C}(A)u_{B,C}-u^2_{B,C}(A) -a_1u_{A,B}+x_{A}(B)+x_{A}(C)+a_2
\end{displaymath}
or, derived from Equation~(\ref{eq:somme}),
\begin{equation}\label{eq:prod}
  u_{A,B}u_{A,C}=x_A+\Gamma(A,B,C)u_{A,C}+\Gamma(A,C,B)u_{A,B} +a_2
+x_{A}(B)+x_{A}(C).
\end{equation}
Indeed, 
\begin{multline*}
(-\frac{1}{z_A}-x_{A}(B)z_A +(y_{A}(B)+a_3)z_A^2)
(-\frac{1}{z_A}-x_{A}(C)z_A +(y_{A}(C)+a_3)z_A^2)+O(z_A^2)\\
=\frac{1}{z_A^2}+x_{A}(B)+x_{A}(C) -(y_{A}(B)+y_{A}(C)+2a_3)z_{A}+O(z_A^2).
\end{multline*}
So, $u_{A,B}u_{A,C}-x_A+a_1u_{A,B}  -x_{A}(B)-x_{A}(C)-a_2$
 cancels at $A$ and its
polar divisor is $-[B]-[C]$. Its residue at $B$ is $-u_{A,B}(C)$.
This proves Equation~(\ref{eq:prod}).
\medskip

In the same vein, we prove 
\begin{equation}\label{eq:carre}
u_{A,B}^2=x_A+x_B-a_1u_{A,B}+x_A(B)+a_2\,.
\end{equation}
Indeed, 
\begin{multline*}
  u_{A,B}^2=(-\frac{1}{z_A}-x_{A}(B)z_A+(y_{A}(B)+a_3)z_A^2)^2+O(z_A^2)\\
  =\frac{1}{z_A^2}+2x_{A}(B)-2(y_{A}(B)+a_3)z_A+O(z_A^2) 
\end{multline*}
and similarly 
\begin{multline*}
 u_{A,B}^2=(\frac{1}{z_B}-a_1+x_{A}(B)z_B+(y_{A}(B)+a_1x_A(B))z_B^2)^2+O(z_B^2)\\
=\frac{1}{z_B^2}-\frac{2a_1}{z_B}+a_1^2+2x_{A}(B)+2y_{A}(B)z_B+O(z_B^2)\,. 
\end{multline*}
So $u_{A,B}^2-x_A-x_B+a_1u_{A,B}=x_{A}(B)+a_2$.
\medskip

Here are more explicit formulas.  
For $A$ and $B$ distinct,
\begin{displaymath}
  u_{A,B} =\left\{
  \begin{array}{ll} 
    -u_{O,A}-a_1 & \text{ if } B = O\,,\\&\\
   {\frac {y+y \left( B \right) +{\it a_1}\,x \left( B \right) +{\it a_3}}{
       x-x \left( B \right) }} & \text{ if } A = O,\\  &\\
    {\frac{{a_1}\,y(A)-3\,x(A)^{2}-2\,{a_2}\,x
       (A)-{a_4}}{2\,y
        (A)+{a_1}\,x(A)+{a_3}}}-{\frac{{a_1}x+{a_3}+2\,y(A)}{x-x(A)}}
    & \text{ if } B = -A\,,\\&\\
   {\frac {y(B)+y(A)+{a_1}\,x(A)+{a_3}}{x(B)-x(A)}}\\&\\
    \hspace*{1cm}
    {+\frac {
        (x(B)-x(A))(y+a_1x+a_3)+
        (y(B)-y(A))x+
        y(A)x(B)-y(B)x(A)
      }
      {(x-x(A))(x-x(B))}
    }
    &\text{ otherwise}.\\   
  \end{array}\right.
\end{displaymath}
Especially, when $A=O$, provided $B$ and $C$ are distinct
and non-zero, we have
\begin{equation}\label{eq:gamma}
\Gamma(O,B,C) = \left\{
  \begin{array}{ll}
    -{\frac {3\,{x(B)}^{2}+{{\it a_1}}\,(y(B)+a_1x(B)+a_3)+2\,{\it a_2}\,x(B)+{\it a_4}}{2\,y(B)+{\it a_1}\,x(B)+{\it a_3}}} & \text{ if } C =
    -B\,,\\ &\\
    \frac{y(C)+y(B)+a_1x(B)+a_3}
    {x(C)-x(B)}
          & \text{ otherwise.}\\ 
  \end{array}\right.
\end{equation}

These formulae can be derived from the definition of $\Gamma  (A,B,C)$
as a slope,  using the explicit form of the addition law on elliptic
curves.

\section{Elliptic bases for finite fields extensions}
\label{section:ellipticbasis}

In this section, we use elliptic functions to construct 
interesting bases for many finite field extensions.
\medskip

Assume  $E$ is an elliptic curve over a finite field
$\bK=\Fq$ and let $d\geqslant 2$ be an integer.
Let $t\in E(\Fq)[d]$ be a rational point of order
$d$. We call $T$ the group generated by $t$. 
Let $\phi: E \rightarrow E$
be the Frobenius endomorphism. 
Let $b\in E(\bar \bK)$  be a point such that $\phi(b)=b+t$.
So $b$   belongs to $E(\bL)$ where $\bL$ is the degree
$d$ extension of $\bK$. 
We denote by $E'$ the quotient $E/T$ and by $I: E\rightarrow E'$
the quotient isogeny. We also assume $db\not =O\in E$.
We set $a=I(b)$ and check $a\in E'(\Fq)$.
For another use of Kummer theory
of  elliptic curves
in order to construct efficient representations for 
finite fields, see \cite{CouLer}.

\subsection{The elliptic basis $\Omega$}

We denote by $\Omega$ the system $(\omega_k)_{k \in \ZZ/d\ZZ}$ defined as
$$\omega_0=1\text{ and }\omega_k=u_{O,kt}(b)\in \bL \text{ for } k\not = 0\bmod d\,.$$
\begin{lemma}
With the above notation, the system 
$\Omega = (\omega_0, \omega_1, \ldots, \omega_{d-1})$ is
 a $\bK$ basis of  $\bL$.
\end{lemma}
\begin{myproof}
Indeed, let the  $\lambda_k$ for $k\in \ZZ/d\ZZ$
be scalars in $\bK$ such that $\sum_{k\in \ZZ/d\ZZ}\lambda_k\omega_k=0$.
The function $f=\lambda_0+\sum_{0\not = k\in \ZZ/d\ZZ}\lambda_ku_{O,kt}$ cancels
at $b$ and also at all its $d$ conjugates over $\bK$ (because $f$ is
defined over $\bK$). But $f$ has no more than $d$ poles (the points
in $T$). If $f$ is non-zero, its  divisor is
 $(f)_0-(f)_\infty$ with $(f)_0=\sum_{t\in T} [b+t]$
and $(f)_\infty=\sum_{t\in T}[t]$. We deduce 
$d\times b$ is zero in $E$. But this is impossible by hypothesis.
Examination of
poles shows that all $\lambda_k$ are zero.
\end{myproof}\\
We call such a basis as $\Omega$ an {\it elliptic basis}.
It enjoys nice properties as we  shall see.
\medskip

We  set
$$\Gamma_{k,l}=\Gamma(O,kt,lt)\in \bK$$
for any distinct non-zero $k, l \in \ZZ/d\ZZ$.
For any  $k\in \ZZ/d\ZZ$, we set  furthermore
$\xi_k=x_{kt}(b)\in \bL$. If $k\not = 0 \bmod d$, we set
 $\nu_k = x_O({kt})\in \bK$ and $\rho_k=y_{O}(kt)\in \bK$ too.

Let now $\Phi: \Fqb\rightarrow \Fqb$ be the $q$-Frobenius
automorphism. We have $x_O(b)=\xi_0$ and
$\Phi(\xi_0)=x_O(\phi(b))=x_O(b+t)=x_{-t}(b)=\xi_{-1}$.
There exist 
$d$ scalars $(\kappa_k)_{0\leqslant k\leqslant d-1}$ in $\bK$ such that 
\begin{equation}\label{eq:reduc}
\xi_0=\sum_{0\leqslant k \leqslant d-1}\kappa_k\omega_k.
\end{equation}
\medskip

We have
for $k\neq 0, 1 \bmod d$,
\begin{eqnarray}\label{eq:conj}
  \Phi(\omega_k)=u_{O,kt}(\phi(b))&=&\nonumber 
  u_{O,kt}(b+t)=u_{-t,(k-1)t}(b)\\&=&\nonumber u_{O,(k-1)t}(b)
  -u_{O,-t}(b)+\Gamma(0,-t,(k-1)t)\\&=&
  \omega_{k-1}-\omega_{-1}+\Gamma_{-1,k-1}
\end{eqnarray}
using Equation~(\ref{eq:somme}).
Similarly 
\begin{equation}\label{eq:conj2}
\Phi(\omega_1)=u_{O, t}(b+t)=u_{-t,O}(b)=-\omega_{-1}-a_1
\text{ and }\Phi(\omega_0)=\omega_0\,.
\end{equation}
Equations~(\ref{eq:conj}) and~(\ref{eq:conj2})
show that the action
of  Frobenius is expressed very easily 
in  an elliptic  basis.
\medskip

As far as multiplication is concerned, we set 
$A=O$, $B=kt$ and $C=lt$
in Equation~(\ref{eq:prod}), and we
evaluate at $b$. We find,
for $k$ and $l$ distinct and non-zero 
in $\ZZ/d\ZZ$,
\begin{equation}\label{eq:mulomega}
  \omega_k\omega_l=\xi_0+\Gamma_{-k,-l}\omega_k+\Gamma_{k,l}\omega_l
  +\nu_k+\nu_{l}+a_2\,.
\end{equation}
In the same vein, from Equation~(\ref{eq:carre}), we obtain for 
any non-zero $k$ in $\ZZ/d\ZZ$, 
\begin{equation}\label{eq:carreomega}
\omega_k^2=\xi_0-a_1\omega_k+\xi_k+\nu_k+a_2\,.
\end{equation}
So, if we multiply two $\bK$-linear combinations of the $\omega$'s, we
quickly get a linear combination of the $\omega$'s and $\xi$'s using
Equations~(\ref{eq:mulomega}) and~(\ref{eq:carreomega}).
We then reduce (eliminate all the $\xi_k$) using the expression of
$\xi_0$ in the basis $\Omega$ given by Equation~(\ref{eq:reduc}). We
also use Equation~(\ref{eq:conj}) to deduce the expressions of all
$\xi_k$'s in the basis $\Omega$.

We don't need to store all constants $\Gamma_{k,l}$.
Equation~(\ref{eq:gamma}) allows to recalculate all these $d^2$ quantities
from the $\nu_k$ and $\rho_k$.  Moreover, we use in the following that only a
small amount of these coefficients has to be computed due to symmetry
relations~(\ref{eq:sym}) and~(\ref{eq:moins}) and invariance by translation.

\paragraph{Example.} Let $\bK=\FF_{7}$ and $d=5$, we first consider the
elliptic curve $E$ of order $10$ defined by
\begin{math}
  {y}^{2}+xy+5\,y={x}^{3}+3\,{x}^{2}+3\,x+2\,.
\end{math}
The point $t=(3,1)$ generates a subgroup $T\subset E$ of order $5$, and with
$E'=E/T$  defined by
\begin{math}
  {y}^{2}+xy+5\,y={x}^{3}+3\,{x}^{2}+4\,x+6\,,
\end{math}
we find
\begin{multline*}
  I: (x,y) \mapsto \left({\frac
      {{x}^{5}+2\,{x}^{2}+5\,x+6}{{x}^{4}+3\,{x}^{2}+4}}, \right.\\
  \left.{\frac { \left( {x}^{6}+4\,{x}^{4}+3\,{x}^{3}+6\,{x}^{2}+3\,x+4 \right) y+3\,{
          x}^{5}+{x}^{4}+{x}^{3}+3\,{x}^{2}+4\,x+1}{{x}^{6}+{x}^{4}+5\,{x}^{2}+6}}\right)\,.
\end{multline*}
Let now $a=(4,2)$, we define $\bL$ with the irreducible polynomial
\begin{math}
  ({\tau}^{5}+2\,{\tau}^{2}+5\,\tau+6)-4\,({\tau}^{4}+3\,{\tau}^{2}+4) = {\tau}^{5}+3\,{\tau}^{4}+4\,{\tau}^{2}+5\,\tau+4\,,
\end{math}
and we set $b=(\tau: \tau^{4756})$.

We find 
\begin{displaymath}
  (u_{O,kt})_{k\in\ZZ/d\ZZ} = \left(1,{\frac {y+2}{x+4}},{\frac {y+2}{x+3}},{\frac {y}{x+3}},{\frac {y+6}{x+4}}\right)\,,  
\end{displaymath}
so that,
\begin{displaymath}
 \Omega=(1, \tau^{10884}, \tau^{11164}, \tau^{9837}, \tau^{15166})\,. 
\end{displaymath}

\subsection{A cell decomposition of the torus}

Equations~(\ref{eq:moins}) and~(\ref{eq:sym}) 
show that the quantity $\Gamma(A,B,C)$ is covariant for the
symmetric group $\cS_3$ and even  for $\cS_3\times \{1,-1\}$.
It is also invariant by translation,
$$\Gamma(A+P,B+P,C+P)=\Gamma(A,B,C).$$
Altogether, $\Gamma$ is covariant for the group
$E(\bar\bK) \rtimes (\cS_3\times \{1,-1\})$.

These covariance properties are useful when
computing the $\Gamma_{k,l}$: we divide by $12$ the amount
of work. 
Since in that case,
 $A=0$, $B=kt$ and $C=lt$ lie  in the group $T=<t>$, a cyclic group
or order $d$,  it makes sense to study the action of $(\ZZ/d\ZZ)\rtimes
(\cS_3\times \{1,-1\})$
on the group $(\ZZ/d\ZZ)^3$. In particular, we are interested in fundamental
domains for this action. It turns out that it is more natural to study first
the action of $\RR^3\rtimes (\cS_3\times \{1,-1\})$ on $\RR^3$.
In this  subsection we justify the choice of fundamental domain that
is made in Subsection~\ref{para:complex}.
\medskip

Let   $\psi: \RR^3\rightarrow \CC$ be the map
that sends the triplet $(a,b,c)$ onto  $a+b\rho+c\rho^2$ where  $\rho=\exp(2i\pi/3)$.
This is a group homomorphism. Its kernel 
is the diagonal subgroup of $\RR^3$. The  group $\cS_3\times \{1, -1\}$  acts
on $\RR^3$ and we have the following covariance formulas
\begin{eqnarray*}
\psi (a,c,b)&=&\overline{\psi(a,b,c)}\,,\\
\psi(c,a,b)&=&\rho\psi(a,b,c)\,,\\
\psi(-a,-b,-c)&=&-\psi(a,b,c)\,.
\end{eqnarray*}
So the map  $\psi$ induces a bijection between the quotient 
of $\RR^3$ by  $\RR \rtimes (\cS_3\times \{1,-1\})$ and the  quotient of 
$\CC$ by  $\mmu_6\times \{1,\conju \}$ where  $\mmu_6$ is 
the group of sixth roots of
unity and  $\conju$ is  complex conjugation.

The image of  $\ZZ^3\subset \RR^3$ by  $\psi$ is the ring of Gaussian integers.
Since  $\ZZ^3$ is normalized by $\cS_3\times \{1,-1\}$, the map 
$\psi$ induces a morphism 
$\tilde \psi: \UU^3 \rightarrow
T_0$ where  $\UU=\RR/\ZZ$ is the unit circle and  $T_0=\CC/(\ZZ+\rho\ZZ)$
the complex torus with zero modular invariant. This map $\tilde \psi$ is
covariant. We denote by $\Lambda$ the lattice $\ZZ+\rho\ZZ$.
For any  $d\geqslant 2$ an integer, we denote by   $\UU[d]$ the  $d$-torsion group
of  $\UU$
and  $T_0[d]$ the one of $T_0$. We denote by  $\psi_d$ the map from $\UU[d]^3$
to  $T_0[d]$ induced by $\tilde \psi$.
\medskip

Let  $k$ and  $l$ be two elements in  $\UU$ and let  $z=k\rho+l\rho^2 \in T_0$
the image of $(0,k,l)$ by $\tilde \psi$. We compute the stabilizer of
$z$  in $\mmu_6\times \{1,\conju \}$. It is clear that $z = \bar z \bmod \Lambda$
  if and only if  $k=l\bmod
1$. The set of fixed points by complex conjugation is the
circle made of real points in $T_0$. 
In the same manner we show that  $-\rho \bar z=   z\bmod \Lambda$
 if and only if $z$ lies
on the circle with equation  $k=2l\bmod 1$. Similarly  
$\rho^2\bar z=z \bmod \Lambda$
if and only if  $l=0\bmod 1$. And  $-\bar z=z \bmod \Lambda$
 if and only if $k=-l\bmod 1$.
And  $\rho\bar z=z \bmod \Lambda $ if and only if 
$k=0\bmod 1$. At last  $-\rho^2 \bar z=z \bmod \Lambda$
 if and only if  $2k=l\bmod 1$.

The only fixed point of $z\bmod \Lambda \mapsto -\rho z \bmod \Lambda$
is  $0$. The same is true for  $z\bmod \Lambda \mapsto -\rho^2z \bmod \Lambda$.

The map  $z\bmod \Lambda \mapsto \rho z\bmod \Lambda $
 has three fixed points, namely $0$,
$(\rho-\rho^2)/{3}$ and its opposite. These are the fixed points of
$z\bmod \Lambda \mapsto \rho^2z\bmod \Lambda$ also. Altogether, these three points form the intersection of the three
circles with equations 
$k=2l\bmod 1$, $l=2k\bmod 1$ and  $l=-k\bmod 1$.

The complementary set of the six  circles above  consists of  $12$ triangles.
Each of these  triangles (with its boundary)
is a fundamental domain for the action of  $\mmu_6\times \{1,\conju\}$
on the torus. The intersection of such a triangle with $T_0[d]$ gives
a fundamental domain for the action of $\mmu_6\times \{1,\conju\}$
on $T_0[d]$. This is also a fundamental domain for the action
of $(\ZZ/d\ZZ)\rtimes (\cS_3\times \{1,-1\})$ on $(\ZZ/d\ZZ)^3$.

\begin{figure}
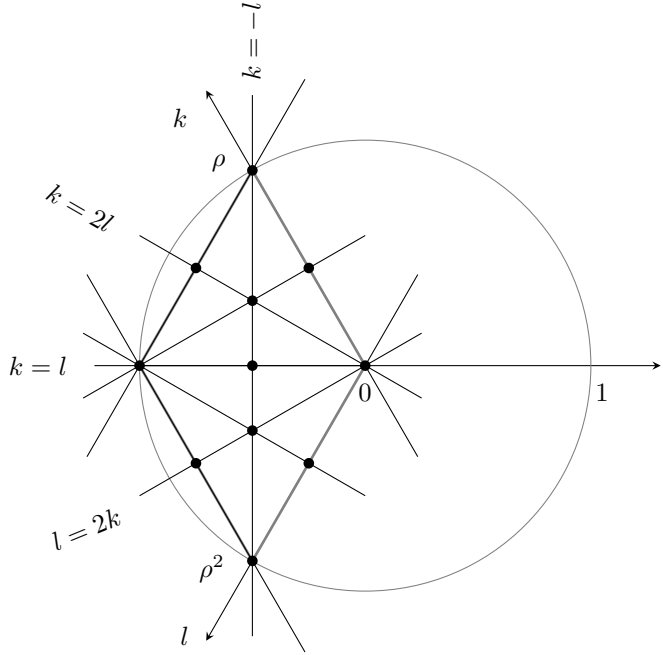

  \centering
  \begin{center}
    \newdimen\gu \gu =3cm
    \begin{pgfpicture}{-2\gu}{-1.3\gu}{2\gu}{1.3\gu}

      \pgfsetlinewidth{0.4pt}
      \pgfsetendarrow{\pgfarrowsingle}
      \pgfline{\pgfpoint{-1.2\gu}{0\gu}}{\pgfpoint{1.3\gu}{0\gu}}
      \pgfputat{\pgfpoint{1.05\gu}{-0.12\gu}}{\pgfbox[center,center]{\footnotesize $1$}}
      \pgfputat{\pgfpoint{-0.0\gu}{-0.12\gu}}{\pgfbox[center,center]{\footnotesize $0$}}

      \pgfline{\pgfpoint{0.233\gu}{-0.404\gu}}{\pgfpoint{-0.7\gu}{1.212\gu}}
      \pgfputat{\pgfpoint{-0.65\gu}{0.9\gu}}{\pgfbox[center,center]{\footnotesize $\rho$}}
      \pgfputat{\pgfpoint{-0.82\gu}{1.1\gu}}{\pgfbox[center,center]{\footnotesize $k$}}

      \pgfline{\pgfpoint{0.233\gu}{0.404\gu}}{\pgfpoint{-0.7\gu}{-1.212\gu}}
      \pgfputat{\pgfpoint{-0.68\gu}{-0.9\gu}}{\pgfbox[center,center]{\footnotesize $\rho^2$}}
      \pgfputat{\pgfpoint{-0.8\gu}{-1.2\gu}}{\pgfbox[center,center]{\footnotesize $l$}}
      \pgfclearendarrow

      {
        \color{DarkGrey}
        \pgfcircle[stroke]{\pgfpoint{0\gu}{0\gu}}{1\gu}
      }

      {
        \color{DarkGrey}
        \pgfsetlinewidth{1pt}
        \pgfline{\pgfpoint{0\gu}{0\gu}}{\pgfpoint{-0.5\gu}{0.866\gu}}
        \pgfline{\pgfpoint{0\gu}{0\gu}}{\pgfpoint{-0.5\gu}{-0.866\gu}}
        \pgfline{\pgfpoint{-0.5\gu}{0.866\gu}}{\pgfpoint{-1\gu}{0\gu}}
        \pgfline{\pgfpoint{-0.5\gu}{-0.866\gu}}{\pgfpoint{-1\gu}{0\gu}}
      }

      {
        \color{black}
        \pgfsetlinewidth{0.4pt}
        \pgfline{\pgfpoint{-0.267\gu}{1.27\gu}}{\pgfpoint{-1.233\gu}{-0.404\gu}}
        \pgfline{\pgfpoint{-0.267\gu}{-1.27\gu}}{\pgfpoint{-1.233\gu}{0.404\gu}}

        \pgfline{\pgfpoint{0\gu}{0\gu}}{\pgfpoint{-1\gu}{0\gu}}
        \pgfputat{\pgfpoint{-1.45\gu}{0pt}}{\pgfbox[center,center]{\footnotesize $k=l$}}

        \pgfline{\pgfpoint{-0.5\gu}{1.2\gu}}{\pgfpoint{-0.5\gu}{-1.2\gu}}
        \pgfputlabelrotated{1.1}%
        {\pgfpoint{-0.5\gu}{-1.2\gu}}{\pgfpoint{-0.5\gu}{1.2\gu}}
        {0pt}{\pgfbox[center,center]{\footnotesize $k=-l$}}

        \pgfline{\pgfpoint{0.25\gu}{0.144\gu}}{\pgfpoint{-1\gu}{-0.577\gu}}
        \pgfputlabelrotated{0.1}%
        {\pgfpoint{-1.5\gu}{-0.577\gu}}{\pgfpoint{0.25\gu}{0.144\gu}}
        {-20pt}{\pgfbox[center,center]{\footnotesize $l=2k$}}

        \pgfline{\pgfpoint{0.25\gu}{-0.144\gu}}{\pgfpoint{-1\gu}{0.577\gu}}
        \pgfputlabelrotated{-0.2}%
        {\pgfpoint{-1\gu}{0.577\gu}}
        {\pgfpoint{0.25\gu}{-0.144\gu}}
        {-2pt}{\pgfbox[center,center]{\footnotesize $k=2l$}}

        \pgfline{\pgfpoint{-1.25\gu}{0.144\gu}}{\pgfpoint{0\gu}{-0.577\gu}}
        \pgfline{\pgfpoint{-1.25\gu}{-0.144\gu}}{\pgfpoint{0\gu}{0.577\gu}}
      }

      {
        \color{black}
        \pgfcircle[fill]{\pgfpoint{0\gu}{0\gu}}{2pt}
        \pgfcircle[fill]{\pgfpoint{-1\gu}{0\gu}}{2pt}

        \pgfcircle[fill]{\pgfpoint{-0.5\gu}{0.866\gu}}{2pt}
        \pgfcircle[fill]{\pgfpoint{-0.5\gu}{-0.866\gu}}{2pt}

        \pgfcircle[fill]{\pgfpoint{-0.25\gu}{0.433\gu}}{2pt}
        \pgfcircle[fill]{\pgfpoint{-0.25\gu}{-0.433\gu}}{2pt}

        \pgfcircle[fill]{\pgfpoint{-0.75\gu}{0.433\gu}}{2pt}
        \pgfcircle[fill]{\pgfpoint{-0.75\gu}{-0.433\gu}}{2pt}

        \pgfcircle[fill]{\pgfpoint{-0.5\gu}{0\gu}}{2pt}
        \pgfcircle[fill]{\pgfpoint{-0.5\gu}{0.288\gu}}{2pt}
        \pgfcircle[fill]{\pgfpoint{-0.5\gu}{-0.288\gu}}{2pt}
      }
    \end{pgfpicture}
  \end{center}
  \caption{Cell decomposition of the torus}
  \label{fig:hexa}
\end{figure}

\subsection{Complexities}\label{para:complex}

Given an elliptic basis $\Omega=(\omega_k)_{k\in\ZZ/d\ZZ}$, we now
focus on the complexity of algorithms for computing the Frobenius or
the multiplication of two elements.  To be as efficient as possible,
and since operands of the algorithms are already of size $d\log q$, we
assume that any precomputation, the storage of which does not exceed
$O(d\log q)$, is possible.\medskip

We first have the following result.
\begin{lemma}
  Let $\alpha =\sum_{i=0}^{d-1} \alpha_i \omega_i \in \bL$. Then there
  exists algorithms that compute $\Phi(\alpha)$ and
  $\Phi^{-1}(\alpha)$ at the expense of $d-1$ multiplications and
  $2d-3$ additions in $\bK$, among which are one multiplication and one
  addition because of the coefficient $a_1$.
\end{lemma}

\begin{figure}[htbp]
  \centering
  \begin{small}
    {\parbox{0.80\textwidth}{\refstepcounter{noalgo}\setcounter{algo}{0}
        \label{algo:EllipticFrobenius}
        \hrule\vskip 0.05cm
        \hrule\vskip 0.1cm
        {\bf Algorithm \thenoalgo} \textsc{EllipticFrobenius}
        \vskip 0.1cm
        \nna Frobenius of an element given in an
        elliptic basis.
        \hrule\vskip 0.1cm
        \nna \textsc{input~:} {$\va = (\alpha_i)_{0\leqslant
            i\leqslant d-1}$ such that $\alpha=\sum_{i=0}^{d-1}
          \alpha_i \omega_i\in \bL$.} 
        \vskip 0cm
        \nna \textsc{output~:} {$\vg=(\gamma_i)_{0\leqslant
            i\leqslant d-1}$ such that $\gamma = \sum_{i=0}^{d-1}
          \gamma_i\omega_i=\Phi(\alpha)\in \bL$.} 
        \hrule\vskip 0.1cm
        \nna \Ret $(\alpha_0-a_1\alpha_1+\sum_{j=2}^{d-1} \alpha_j
        \Gamma_{d-1,j-1}, \alpha_2, \ldots, \alpha_{d-1},
        -\sum_{j=1}^{d-1} \alpha_j )$
        \hrule\vskip 0.5cm
      }
    }
  \end{small}  
\end{figure}

\begin{figure}[htbp]
  \centering
  \begin{small}
    {\parbox{0.80\textwidth}{\refstepcounter{noalgo}\setcounter{algo}{0}
        \label{algo:EllipticFrobeniusInverse}
        \hrule\vskip 0.05cm
        \hrule\vskip 0.1cm
        {\bf Algorithm \thenoalgo} \textsc{EllipticFrobeniusInverse}
        \vskip 0.1cm
        \nna Inverse Frobenius of an element given in an elliptic basis.
        \hrule\vskip 0.1cm
        \nna \textsc{input~:} {$\va = (\alpha_i)_{0\leqslant
            i\leqslant d-1}$ such that $\alpha=\sum_{i=0}^{d-1}
          \alpha_i \omega_i\in \bL$.}
        \vskip 0cm
        \nna \textsc{output~:} {$\vg=(\gamma_i)_{0\leqslant
            i\leqslant d-1}$ such that $\gamma = \sum_{i=0}^{d-1}
          \gamma_i\omega_i=\Phi^{-1}(\alpha)\in \bL$.}
        \hrule\vskip 0.1cm
        \nna \Ret $
        (\alpha_0+\sum_{j=1}^{d-2} \alpha_j \Gamma_{j,d-1}-a_1\alpha_{d-1},
        -\sum_{j=1}^{d-1} \alpha_j,
        \alpha_1, \ldots, \alpha_{d-2})
        $ \hrule\vskip 0.5cm } }
  \end{small}  
\end{figure}

\begin{myproof} Plugging Equation~(\ref{eq:conj}) and
  Equation~(\ref{eq:conj2}) in $\sum_{i=0}^{d-1} \alpha_i \Phi(\omega_i)$ or
  $\sum_{i=0}^{d-1} \alpha_i \Phi^{-1}(\omega_i)$ proves the correctness of
  Algorithm~\ref{algo:EllipticFrobenius} and
  Algorithm~\ref{algo:EllipticFrobeniusInverse}.  And, once precomputed the
  $\Gamma_{d-1,j}$'s and $\Gamma_{j,d-1}$'s, the complexity is obvious.
\end{myproof}

Multiplying two elements in such a basis can be done with good complexity too.
\begin{lemma}\label{le:mult}
  Let $\alpha =\sum_{i=0}^{d-1} \alpha_i \omega_i \in \bL$ and $\beta
  =\sum_{i=0}^{d-1} \beta_i \omega_i \in \bL$. Then there exists an
  algorithm that computes the product $\alpha \times \beta$ at the expense of
  \begin{itemize}
  \item $(37\,d^2+30\,d-7\varepsilon-60)/12$ additions,
    $(32\,d^2+42\,d-2\varepsilon-48)/12$ multiplications and
    $(d^2-\varepsilon)/12$ inversions in $\bK$,
  \end{itemize}
  where $\varepsilon=12,1,4,9,4,1$ respectively for $d=0,\ldots,5\bmod 6$,
  among which are
  $(d^2+12d-\varepsilon-24)/12$ additions and
    $(d^2+36\,d-\varepsilon-48)/12$ multiplications because of the coefficient $a_1$,
  $(d^2-\varepsilon)/12$ additions because of the coefficient $a_3$.
\end{lemma}

\begin{figure}[htbp]
  \centering
\begin{small}
{\parbox{\textwidth}{\refstepcounter{noalgo}\setcounter{algo}{0}
  \label{algo:EllipticMultiplication}
  \hrule\vskip 0.05cm
  \hrule\vskip 0.1cm
  {\bf Algorithm \thenoalgo} \textsc{EllipticMultiplication}
  \vskip 0.1cm
  \nna Product of two elements given in an elliptic basis.
  \hrule\vskip 0.1cm
  \nna \textsc{input~:} {$\va = (\alpha_i)_{0\leqslant
      i\leqslant d-1}$ and $\vb = (\beta_i)_{0\leqslant
      i\leqslant d-1}$ such that $\alpha=\sum_{i=0}^{d-1}
    \alpha_i \omega_i,$ $\beta=\sum_{i=0}^{d-1}
    \beta_i \omega_i\in \bL$.} 
  \vskip 0cm
        \nna \textsc{output~:} {$\vg=(\gamma_i)_{0\leqslant
            i\leqslant d-1}$ such that $\gamma = \sum_{i=0}^{d-1}
          \gamma_i\omega_i=\alpha \times \beta\in \bL$.} 
  \hrule\vskip 0.1cm
  \nns{} $s_a:= 0$\,; $s_b:= \beta_1$\,; $\gamma_0:= 0$\,; $\gamma_1:= -a_1 s_b \alpha_1$\,;
  \nns{} \For $k:= 2$ \To $d-1$ \Do
  $s_a +\!\!:= \alpha_{k-1}$\,; $s_b +\!\!:= \beta_k$\,; $\gamma_k:= -a_1 (s_b \alpha_k+s_a\beta_k)$\,;
  \Endindent{}
  \nns{} $s_a +\!\!:= \alpha_{d-1}$\,; $(\gamma_0,\ldots, \gamma_{d-1}) +\!\!:=
  s_as_b\,(\kappa_0+a_2,\kappa_1, \ldots \kappa_{d-1})$\,;
  \nns{} $s_a':= \sum_{i=1}^{d-1} \alpha_i\nu_i$\,; $s_b':= \sum_{i=1}^{d-1}
  \beta_i\nu_i$\,; $\gamma_0~+\!\!:= s_as_b'+s_a's_b$\,;
  \nns{} \For $k:= 1$ \To $d-1$ \Do
  \nns{}  $\delta:= \alpha_k\beta_k$\,;
  $\gamma_{0} +\!\!:= \delta\,((\Phi^{-k}(\xi_0))_0 - \nu_k)$\,;
  $\gamma_{k} -\!\!:= \delta\,\sum_{l=1}^{d-1} \kappa_l$;
  \nns{} \For $l:= 1$ \To $k-1$ \Do $\gamma_{l} +\!\!:= \delta\,\kappa_{(d-k+l) \bmod d}$;\Endindent{}
  \nns{} \For $l:= k+1$ \To $d-1$ \Do $\gamma_{l} +\!\!:= \delta\,\kappa_{(d-k+l) \bmod d}$;\Endindent{}
  \Endindent{}
  \nns{} $(\gamma_0,\ldots, \gamma_{d-1})~+\!\!:= (\alpha_0\beta_0,
  \alpha_1\beta_0+\alpha_0\beta_1,\ldots,\alpha_{d-1}\beta_0+\alpha_0\beta_{d-1})$\,;
  \nns{} \If $d \bmod 3 = 0$ \Then
  \nns{} $g:= - (3\,\nu_{2d/3}^2 + 2 a_2 \nu_{2d/3} + a_4) / (2\rho_{2d/3}+a_1\nu_{2d/3}+a_3)-a_1$\,; 
  \nns{} $\delta:= g\,(\alpha_{2d/3}\beta_{d/3} +
  \alpha_{d/3}\beta_{2d/3})$\,; $\gamma_{2d/3}~-\!\!:= \delta$\,; $\gamma_{d/3}~+\!\!:= \delta$\,;
  \Endindent{}
  \nns{} \For $k:= 2$ \To $\lfloor (2d-1)/3\rfloor$ \By $2$ \Do
  \nns{} $l:= k/2$\,; $g:= (\rho_l+\rho_k+a_1\nu_k+a_3)/(\nu_l-\nu_k)$\,;
  \nns{} $i_1,i_2:= 2\,l,d-l$\,; $j_1,j_2:= d-2\,l, l$\,;
  \nns{} $\delta_{12}:= g\,(\alpha_{i_1}\,\beta_{j_2}+\alpha_{j_2}\,\beta_{i_1})$\,;
  $\delta_{21}:= g\,(\alpha_{i_2}\,\beta_{j_1}+\alpha_{j_1}\,\beta_{i_2})$\,;
  $\delta_{22}:= g\,(\alpha_{i_2}\,\beta_{j_2}+\alpha_{j_2}\,\beta_{i_2})$\,;
  \nns{} $\gamma_{i_1}~-\!\!:= \delta_{12}$\,;
  $\gamma_{i_2}~-\!\!:= \delta_{21} + \delta_{22}$\,;
  $\gamma_{j_1}~+\!\!:= \delta_{21}$\,; 
  $\gamma_{j_2}~+\!\!:= \delta_{12} + \delta_{22}$\,;
  \Endindent{}
  \nns{} \For $k:= \lfloor 1+d/2\rfloor$ \To  $\lfloor (2d-1)/3\rfloor$ \Do
  \nns{}  $l:= 2k \bmod d$\,; $g:= (\rho_l+\rho_k+a_1\nu_k+a_3)/(\nu_l-\nu_k)$\,;
  \nns{}  $i_1, i_2:=  k, (2d - 2k) \bmod d$\,; $j_1,j_2:= (2k) \bmod d, d - k$;
  \nns{} $\delta_{11}:= g\,(\alpha_{i_1}\,\beta_{j_1}+\alpha_{j_1}\,\beta_{i_1})$\,;
  $\delta_{22}:= g\,(\alpha_{i_2}\,\beta_{j_2}+\alpha_{j_2}\,\beta_{i_2})$\,;
  $\delta_{12}:= g\,(\alpha_{i_1}\,\beta_{j_2}+\alpha_{j_2}\,\beta_{i_1})$\,;
  \nns{} $\gamma_{i_1}~-\!\!:= \delta_{11} + \delta_{12}$\,;
  $\gamma_{i_2}~-\!\!:= \delta_{22}$\,;
  $\gamma_{j_1}~+\!\!:= \delta_{11}$\,;
  $\gamma_{j_2}~+\!\!:= \delta_{22} + \delta_{12}$\,;
  \Endindent{}
  \nns{} \For $k:= 3$ \To $\lfloor (2d-1)/3\rfloor$ \Do
  \nns{} \For $l:= \max(1, 2k-d + 1)$ \To $\lfloor (k-1)/2\rfloor$ \Do
  \nns{}  $g:= (\rho_l+\rho_k+a_1\nu_k+a_3)/(\nu_l-\nu_k)$\,;
  \nns{}  $i_1, i_2, i_3:= k, d-l, d-k+l$\,; $j_1, j_2, j_3:= d-k, l, k-l$;
  \nns{}
	 $\delta_{12}:= g\,(\alpha_{i_1}\,\beta_{j_2} + \alpha_{j_2}\,\beta_{i_1})$\,;
	 $\delta_{13}:= g\,(\alpha_{i_1}\,\beta_{j_3} + \alpha_{j_3}\,\beta_{i_1})$\,;
	 $\delta_{21}:= g\,(\alpha_{i_2}\,\beta_{j_1} + \alpha_{j_1}\,\beta_{i_2})$\,;
  \nns{} 
	 $\delta_{23}:= g\,(\alpha_{i_2}\,\beta_{j_3} + \alpha_{j_3}\,\beta_{i_2})$\,;
	 $\delta_{31}:= g\,(\alpha_{i_3}\,\beta_{j_1} + \alpha_{j_1}\,\beta_{i_3})$\,;
	 $\delta_{32}:= g\,(\alpha_{i_3}\,\beta_{j_2} + \alpha_{j_2}\,\beta_{i_3})$\,;
  \nns{}
	 $\gamma_{i_1}~-\!\!:= \delta_{12} + \delta_{13}$\,;
         $\gamma_{i_2}~-\!\!:= \delta_{21} + \delta_{23}$\,;
         $\gamma_{i_3}~-\!\!:= \delta_{31} + \delta_{32}$\,;
  \nns{}
	 $\gamma_{j_1}~+\!\!:= \delta_{21} + \delta_{31}$\,;
	 $\gamma_{j_2}~+\!\!:= \delta_{12} + \delta_{32}$\,;
         $\gamma_{j_3}~+\!\!:= \delta_{13} + \delta_{23}$\,;
  \Endindent{}
  \Endindent{}
  \nns{} \Ret $(\gamma_i)_{0\leqslant i\leqslant d-1}$
  \vskip 0.1cm\hrule\vskip 0.5cm
}}
\end{small}
\end{figure}

\begin{myproof} We prove the correctness of
  Algorithm~\ref{algo:EllipticMultiplication} and establish its
  complexity.\medskip

  \noindent\hskip0.5cm\textit{Correctness.}
  Equations~(\ref{eq:prod}) and~(\ref{eq:carre}), for $k \leqslant l$, yield
  \begin{displaymath}
    \omega_k\, \omega_l = \omega_l\, \omega_k = \left\{
      \begin{array}{ll}
        \omega_l & \text{if }k = 0\,,\\
        \xi_0+a_2-a_1\omega_k+\Phi^{-k}(\xi_0)+\nu_k\,\omega_0 & \text{if } l = k\text{ and }k > 0\,,\\
        \xi_0+a_2-a_1\omega_k+\Gamma_{k, l}\,(\omega_l-\omega_k)+(\nu_k+\nu_l)\,\omega_0 &
        \text{otherwise\,.}
      \end{array}\right.
  \end{displaymath}
  And we have,
  \begin{multline}\label{eq:sums}
    \alpha \times \beta = 
    \sum_{k=0}^{d-1} \sum_{l=0}^{d-1} \alpha_k\beta_l \omega_k \omega_l = 
    (\sum_{k=1}^{d-1} \alpha_k)(\sum_{l=1}^{d-1} \beta_l)(\xi_0+a_2)\\
    +\left((\sum_{k=1}^{d-1} \alpha_k)(\sum_{l=1}^{d-1} \beta_l\nu_l)
    +(\sum_{k=1}^{d-1} \alpha_k\nu_k)(\sum_{l=1}^{d-1} \beta_l)\right)\omega_0\\
    +\alpha_0\beta_0\omega_0+\sum_{k=1}^{d-1} \alpha_k\,\beta_k\,(\Phi^{-k}(\xi_0)-\nu_k\,\omega_0)
    +\sum_{k=1}^{d-1} (\alpha_k\beta_0+\beta_k\alpha_0)\omega_k\\ 
    -a_1\sum_{0 < k, l <d} \alpha_k\beta_l \omega_k
    +\sum_{\stackrel{0 < k, l < d}{k \neq l}}
    \Gamma_{k,l}\alpha_k\beta_l(\omega_l-\omega_k)\,.
  \end{multline}
  The first two terms of this sum are computed at steps 3. and 4. of the
  algorithm. The three next terms are computed in steps 5. to 9. Especially,
  steps 5. to 8. correspond to the action of $\Phi^{-k}$ on $\xi_0$ (the
  quantity $(\Phi^{-k}(\xi_0))_0$, at step 4., is the first coordinate of
  $\Phi^{-k}(\xi_0)$ written in basis $\Omega$). 

  The constants $\Gamma_{k,l}$ satisfied 12 symmetry relations and we take
  advantage of them to compute the two last terms of the sum.  More precisely,
  for $k$ and $l$ distinct and non-zero in $\ZZ/d\ZZ$, we have
  \begin{displaymath}
    \left\{
    \begin{array}{l}
      \Gamma_{k,l} = \Gamma_{-l,-k} = \Gamma_{k,k-l} = \Gamma_{l-k,-k} =
      \Gamma_{l-k,l} = \Gamma_{-l,k-l}\,,\\
      \Gamma_{l,k} = \Gamma_{-k,-l} = \Gamma_{k-l,k} = \Gamma_{-k,l-k} = \Gamma_{l,l-k} =
      \Gamma_{k-l,-l}\,, 
    \end{array}
    \right.\text{ and } \Gamma_{k,l} = -\Gamma_{l,k} - a_1\,. 
  \end{displaymath}
  All of these relations can be proved thanks to Equation~(\ref{eq:sym})  and
  Equation~(\ref{eq:moins}). For instance,
  to check that $\Gamma_{k,l} = \Gamma_{l-k,-k}$, we start from
  \begin{math}
    \Gamma(O,kt,lt) = u_{O,kt}(b+kt)+u_{kt,lt}(b+kt)+u_{lt,O}(b+kt),
  \end{math}
  and we find 
  \begin{math}
    \Gamma(O,kt,lt) = u_{-kt,O}(b)+u_{O,(l-k)t}(b)+u_{(l-k)t,-kt}(b)=\Gamma(O,(l-k)t,-kt)\,.
  \end{math}

  \begin{figure}[htb]
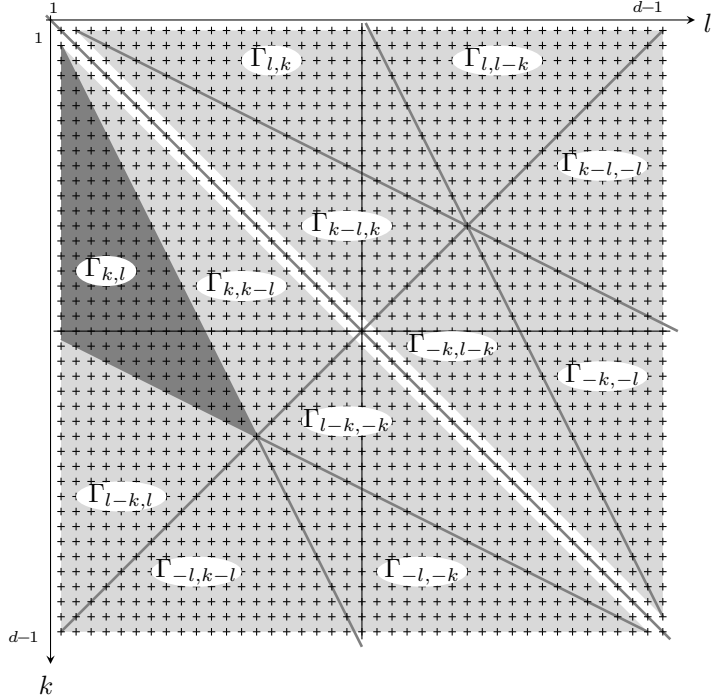

    \begin{center}
      \newdimen\gu   \gu =0.20cm
      \newdimen\tcrux
      \newcount\cntk \newcount\cntl \

      \begin{pgfpicture}{-2\gu}{-2\gu}{44\gu}{44\gu}

        {\color{DarkGrey}
          \pgfmoveto{\pgfpoint{1\gu}{40\gu}}
          \pgflineto{\pgfpoint{1\gu}{20.5\gu}}
          \pgflineto{\pgfpoint{14\gu}{14\gu}}
          \pgffill
        }
        {\color{LightGrey}
          \pgfmoveto{\pgfpoint{1\gu}{20.5\gu}}
          \pgflineto{\pgfpoint{14\gu}{14\gu}}
          \pgflineto{\pgfpoint{1\gu}{1\gu}}
          \pgffill
        }
        {\color{LightGrey}
          \pgfmoveto{\pgfpoint{1\gu}{1\gu}}
          \pgflineto{\pgfpoint{14\gu}{14\gu}}
          \pgflineto{\pgfpoint{20.5\gu}{1\gu}}
          \pgffill
        }
        {\color{LightGrey}
          \pgfmoveto{\pgfpoint{20.5\gu}{1\gu}}
          \pgflineto{\pgfpoint{14\gu}{14\gu}}
          \pgflineto{\pgfpoint{40\gu}{1\gu}}
          \pgffill
        }
        {\color{LightGrey}
          \pgfmoveto{\pgfpoint{40\gu}{1\gu}}
          \pgflineto{\pgfpoint{14\gu}{14\gu}}
          \pgflineto{\pgfpoint{20.5\gu}{20.5\gu}}
          \pgffill
        }
        {\color{LightGrey}
          \pgfmoveto{\pgfpoint{14\gu}{14\gu}}
          \pgflineto{\pgfpoint{20.5\gu}{20.5\gu}}
          \pgflineto{\pgfpoint{1\gu}{40\gu}}
          \pgffill
        }

        {\color{LightGrey}
          \pgfmoveto{\pgfpoint{2\gu}{41\gu}}
          \pgflineto{\pgfpoint{21.5\gu}{41\gu}}
          \pgflineto{\pgfpoint{28\gu}{28\gu}}
          \pgffill
        }
        {\color{LightGrey}
          \pgfmoveto{\pgfpoint{21.5\gu}{41\gu}}
          \pgflineto{\pgfpoint{28\gu}{28\gu}}
          \pgflineto{\pgfpoint{41\gu}{41\gu}}
          \pgffill
        }
        {\color{LightGrey}
          \pgfmoveto{\pgfpoint{41\gu}{41\gu}}
          \pgflineto{\pgfpoint{28\gu}{28\gu}}
          \pgflineto{\pgfpoint{41\gu}{21.5\gu}}
          \pgffill
        }
        {\color{LightGrey}
          \pgfmoveto{\pgfpoint{28\gu}{28\gu}}
          \pgflineto{\pgfpoint{41\gu}{21.5\gu}}
          \pgflineto{\pgfpoint{41\gu}{2\gu}}
          \pgffill
        }
        {\color{LightGrey}
          \pgfmoveto{\pgfpoint{41\gu}{2\gu}}
          \pgflineto{\pgfpoint{28\gu}{28\gu}}
          \pgflineto{\pgfpoint{21.5\gu}{21.5\gu}}
          \pgffill
        }
        {\color{LightGrey}
          \pgfmoveto{\pgfpoint{21.5\gu}{21.5\gu}}
          \pgflineto{\pgfpoint{28\gu}{28\gu}}
          \pgflineto{\pgfpoint{2\gu}{41\gu}}
          \pgffill
        }

        \pgfsetlinewidth{1pt}

        {
          \color{DarkGrey}
          \pgfline{\pgfpoint{1\gu}{40\gu}}{\pgfpoint{21\gu}{0\gu}}
          \pgfline{\pgfpoint{1\gu}{20.5\gu}}{\pgfpoint{40\gu}{1\gu}}
        }

        {
          \color{DarkGrey}
          \pgfline{\pgfpoint{2\gu}{41\gu}}{\pgfpoint{42\gu}{21\gu}}
          \pgfline{\pgfpoint{21.25\gu}{41.5\gu}}{\pgfpoint{41\gu}{2\gu}}
        }

        {\color{DarkGrey}
          \pgfline{\pgfpoint{41\gu}{41\gu}}{\pgfpoint{1\gu}{1\gu}}
        }

        {\color{DarkGrey}
          \pgfline{\pgfpoint{0\gu}{42\gu}}{\pgfpoint{41.5\gu}{0.5\gu}}
        }

        \pgfsetlinewidth{0.4pt}

        \cntk=1
        \loop\ifnum\cntk<42 {
          \cntl=1
          \loop\ifnum\cntl<42

             \tcrux=\cntk\gu \advance \tcrux by 0.25\gu
             \pgfline{\pgfpoint{\cntk\gu}{\cntl\gu}}{\pgfpoint{\tcrux}{\cntl\gu}}
             \tcrux=\cntk\gu \advance \tcrux by -0.25\gu
             \pgfline{\pgfpoint{\cntk\gu}{\cntl\gu}}{\pgfpoint{\tcrux}{\cntl\gu}}

             \tcrux=\cntl\gu \advance \tcrux by 0.25\gu
             \pgfline{\pgfpoint{\cntk\gu}{\cntl\gu}}{\pgfpoint{\cntk\gu}{\tcrux}}
             \tcrux=\cntl\gu \advance \tcrux by -0.25\gu
             \pgfline{\pgfpoint{\cntk\gu}{\cntl\gu}}{\pgfpoint{\cntk\gu}{\tcrux}}

          \advance \cntl by 1
          \repeat
        }
        \advance \cntk by 1
        \repeat

        \pgfsetendarrow{\pgfarrowsingle}
        \pgfline{\pgfpoint{0.3\gu}{42\gu}}{\pgfpoint{0.3\gu}{-1\gu}}
        \pgfputat{\pgfpoint{0\gu}{-2.5\gu}}{\pgfbox[center,center]{\footnotesize $k$}}
        \pgfline{\pgfpoint{0\gu}{41.7\gu}}{\pgfpoint{43\gu}{41.7\gu}}
        \pgfputat{\pgfpoint{44\gu}{41.5\gu}}{\pgfbox[center,center]{\footnotesize $l$}}
        \pgfclearendarrow

        \pgfputat{\pgfpoint{0.5\gu}{42.5\gu}}{\pgfbox[center,center]{\tiny $1$}}
        \pgfputat{\pgfpoint{40\gu}{42.5\gu}}{\pgfbox[center,center]{\tiny $d\!\!-\!\!1$}}
        \pgfputat{\pgfpoint{-0.5\gu}{40.5\gu}}{\pgfbox[center,center]{\tiny $1$}}
        \pgfputat{\pgfpoint{-1.5\gu}{0.6\gu}}{\pgfbox[center,center]{\tiny $d\!\!-\!\!1$}}

        \pgfline{\pgfpoint{0.5\gu}{21\gu}}{\pgfpoint{41.5\gu}{21\gu}}
        \pgfline{\pgfpoint{21\gu}{41.5\gu}}{\pgfpoint{21\gu}{0.5\gu}}

        {\color{white}
          \pgfellipse[fill]{\pgfpoint{4\gu}{25\gu}}{\pgfxy(0.4,0)}{\pgfxy(0,0.2)}
          \pgfellipse[fill]{\pgfpoint{15\gu}{39\gu}}{\pgfxy(0.4,0)}{\pgfxy(0,0.2)}
          \pgfellipse[fill]{\pgfpoint{25\gu}{5\gu}}{\pgfxy(0.6,0)}{\pgfxy(0,0.2)}
          \pgfellipse[fill]{\pgfpoint{37\gu}{18\gu}}{\pgfxy(0.6,0)}{\pgfxy(0,0.2)}
          \pgfellipse[fill]{\pgfpoint{27\gu}{20\gu}}{\pgfxy(0.6,0)}{\pgfxy(0,0.2)}
          \pgfellipse[fill]{\pgfpoint{20\gu}{15\gu}}{\pgfxy(0.6,0)}{\pgfxy(0,0.2)}
          \pgfellipse[fill]{\pgfpoint{20\gu}{28\gu}}{\pgfxy(0.6,0)}{\pgfxy(0,0.2)}
          \pgfellipse[fill]{\pgfpoint{13\gu}{24\gu}}{\pgfxy(0.6,0)}{\pgfxy(0,0.2)}
          \pgfellipse[fill]{\pgfpoint{30\gu}{39\gu}}{\pgfxy(0.6,0)}{\pgfxy(0,0.2)}
          \pgfellipse[fill]{\pgfpoint{5\gu}{10\gu}}{\pgfxy(0.6,0)}{\pgfxy(0,0.2)}
          \pgfellipse[fill]{\pgfpoint{37\gu}{32\gu}}{\pgfxy(0.6,0)}{\pgfxy(0,0.2)}
          \pgfellipse[fill]{\pgfpoint{10\gu}{5\gu}}{\pgfxy(0.6,0)}{\pgfxy(0,0.2)}
        }

        {\color{black}
          \pgfputat{\pgfpoint{4\gu}{25\gu}}{\pgfbox[center,center]{\footnotesize $\Gamma_{k,l}$}}
          \pgfputat{\pgfpoint{15\gu}{39\gu}}{\pgfbox[center,center]{\footnotesize $\Gamma_{l,k}$}}
          \pgfputat{\pgfpoint{25\gu}{5\gu}}{\pgfbox[center,center]{\footnotesize $\Gamma_{-l,-k}$}}
          \pgfputat{\pgfpoint{37\gu}{18\gu}}{\pgfbox[center,center]{\footnotesize $\Gamma_{-k,-l}$}}
          \pgfputat{\pgfpoint{27\gu}{20\gu}}{\pgfbox[center,center]{\footnotesize $\Gamma_{-k,l-k}$}}
          \pgfputat{\pgfpoint{20\gu}{15\gu}}{\pgfbox[center,center]{\footnotesize $\Gamma_{l-k,-k}$}}
          \pgfputat{\pgfpoint{20\gu}{28\gu}}{\pgfbox[center,center]{\footnotesize $\Gamma_{k-l,k}$}}
          \pgfputat{\pgfpoint{13\gu}{24\gu}}{\pgfbox[center,center]{\footnotesize $\Gamma_{k,k-l}$}}
          \pgfputat{\pgfpoint{30\gu}{39\gu}}{\pgfbox[center,center]{\footnotesize $\Gamma_{l,l-k}$}}
          \pgfputat{\pgfpoint{5\gu}{10\gu}}{\pgfbox[center,center]{\footnotesize $\Gamma_{l-k,l}$}}
          \pgfputat{\pgfpoint{37\gu}{32\gu}}{\pgfbox[center,center]{\footnotesize $\Gamma_{k-l,-l}$}}
          \pgfputat{\pgfpoint{10\gu}{5\gu}}{\pgfbox[center,center]{\footnotesize $\Gamma_{-l,k-l}$}}
        }
      \end{pgfpicture}
    \end{center}
    \caption{Symmetry relations on the coefficients $\Gamma_{k,l}$ {\small ($d=42$)}}
    \label{fig:gammadomain}
  \end{figure}

  We use first that $\Gamma_{k,l} = -\Gamma_{l,k} - a_1$ and we rewrite the last
  two terms of Equation~(\ref{eq:sums}) as follows,
  \begin{displaymath}
    -a_1\sum_{k=1}^{d-1} (\beta_k\sum_{l=1}^{k-1}\alpha_l+ \alpha_k\sum_{l=1}^{k}\beta_l) \omega_k
    +\sum_{0 < l < k < d}
    \Gamma_{k,l}(\alpha_k\beta_l+\alpha_l\beta_k)(\omega_l-\omega_k)\,.
  \end{displaymath}
  The first term of this sum is computed at at steps 1. and 2. of the
  algorithm.  To compute the last term, we consider in turn each orbit
  of the action defined by the symmetries on the coefficients
  $\Gamma_{k,l}$.  We choose as a fundamental domain for this action
  the triangle delimited by the circles $l=1$, $k = 2\,l \bmod d$ and
  $l = 2\,k \mod d$ (cf.  Figure~\ref{fig:gammadomain}). It is
  cumbersome, but not difficult, to check that any point of this
  domain, outside the two circles $k = 2\,l \bmod d$ and $l = 2\,k
  \mod d$, has an orbit of exactly $12$ points: we compute only once the
  constant $\Gamma_{k,l}$ corresponding to these $12$ points and we
  calculate accordingly their contribution to the product
  $\alpha\times \beta$. These are steps 23. to 30. of the algorithm.

  Points on the line $k = 2\,l \bmod d$ have orbits of only $6$ points. We
  precisely have
  \begin{math}
    \Gamma_{2\,l, l} = \Gamma_{-l, -2\,l} = \Gamma_{-l, l} =
    -\Gamma_{l, 2\,l} -a_1= -\Gamma_{-2\,l, -l} -a_1= -\Gamma_{l, -l}-a_1\,,
  \end{math}
  and this yield steps 13. to 17. of the algorithm.  Similarly, points on the
  line $l = 2\,k \bmod d$ have orbits of only $6$ points too. We have
  \begin{math}
    \Gamma_{k, 2\,k} = \Gamma_{-2\,k,-k} = \Gamma_{k, -k} =
    -\Gamma_{2\,k, k} -a_1= -\Gamma_{-k, -2\,k} -a_1= -\Gamma_{-k, k} -a_1
  \end{math}
  and this yield steps 18. to 22. of the algorithm.
  
  Finally, when $d$ is divisible by 3, the two circles $k = 2\,l \bmod d$ and $l =
  2\,k \bmod d$ meet at the exceptional point $(2d/3, d/3)$, which is on the
  $k+l=0\bmod d$ line too. This point has an orbit of only 2 points,
  \textit{i.e.}
  \begin{math}
    \Gamma_{2d/3, d/3} = -\Gamma_{d/3, 2d/3}-a_1\,.
  \end{math}
  This yields steps 10. to 12. of the algorithm.
  \medskip

  \noindent\hskip0.5cm\textit{Complexity.} We precompute the
  $d$ constants $\nu_k$ and $\rho_k$, the constant $\Gamma_{2d/3, d/3}$ if
  $d\bmod 3 = 0$, the $d$ coordinates in the basis $\Omega$ of $\xi_0$, their
  sum $\sum_{l=1}^{d-1} \kappa_l$, $\kappa_0+a_2$ and the
  $\omega_0$-coordinates of all $\Phi^k(\xi_0)-\nu_k$ for $0\leqslant k\leqslant d-1$.

  Then,
  Steps 1.-2. need $3d-7$ additions and $3d-4$ multiplications in $\bK$ (among
  which are
  $d-2$ additions and $3d-4$ multiplications because of $a_1$),
  Step 3. needs $d+1$ additions and $d+1$ multiplications in $\bK$,
  Step 6. needs $d-1$ additions and $2d-2$ multiplications in $\bK$,
  Steps 7.-8. need $d^2-2d+1$ additions and $d^2-2d+1$ multiplications in $\bK$,
  Step 9. needs $2d-1$ additions and $2d-1$ multiplications in $\bK$,
  Steps 11.-12. need $3$ additions and $3$ multiplications in $\bK$ if $d$ is a
  multiple of 3 (and cost nothing otherwise),
  Steps 13.-17. consist in $\lfloor(d-1)/3\rfloor$ iterations and Steps
  18.-22. consist in $\lfloor(d-5+6\varepsilon')/6\rfloor$ (where $\varepsilon'
  = 0$ if $d \bmod 6 = 0$ and $\varepsilon' = 1$ otherwise), each of them needs
  $16$ additions, $11$ multiplications and $1$ inversion in $\bK$ (among which
  are $1$
  addition, $1$ multiplication because of $a_1$ and $1$ addition because of $a_3$),
  and finally, Steps 23.-30. consist in $\lfloor d^2/12 \rfloor -\lfloor
  d/2\rfloor+\varepsilon''$ iterations (where $\varepsilon'' = 0$ if $d
  \bmod 6 = 1,5$ and $\varepsilon'' = 1$ otherwise), each of them needs
  $25$ additions, $12$ multiplications and $1$ inversion in $\bK$ (among which
  are $1$
  addition, $1$ multiplication because of $a_1$ and $1$ addition because of $a_3$).

  Adding all these complexities yields the complexity announced.
\end{myproof}

Depending on the characteristic of $\bK$, it is classical to consider the reduced
Weierstrass Model to define elliptic curves. We give in Table~\ref{tab:ec}
precise complexities for these cases, all obtained with Lemma~\ref{le:mult}.
\begin{figure}[htbp]
  \begin{small}
    \begin{displaymath}%
      \begin{array}{|ll||c|c|c|c|}
        \hline
        \multicolumn{2}{|c||}{\text{Condition}} & \text{Model} & \text{Add.} &
        \text{Mult.} & \text{Inv.}\\
        \hline \hline
        \Car(\bK)\neq 2,3 & &Y^2=X^3+a_4X+a_6 &&&\\
        \Car(\bK) = 3, & j_E \neq 0 & Y^2=X^3+a_2X^2+a_6 &
        \frac{35\,{d}^{2}+18\,d-5\,\varepsilon-36}{12} &
        &\\
        & j_E = 0 & Y^2=X^3+a_4X+a_6 &&
        \frac{31\,{d}^{2}+6\,d-\varepsilon}{12} &
        \frac{{d}^{2}-\varepsilon}{12}\\
        \cline{1-4}
        \Car(\bK) = 2 & j_E \neq 0 & Y^2+XY=X^3+a_2X^2+a_6 &
        \frac{6\,{d}^{2}+5\,d-\varepsilon-10}{2} &&\\
        & j_E = 0 & Y^2+a_3Y=X^3+a_4X+a_6 &
        \frac{6\,{d}^{2}+3\,d-\varepsilon-6}{2}&&\\
        \hline
      \end{array}    
    \end{displaymath}
  \end{small}
  \begin{center}
    \caption{Elliptic multiplication complexities}
    \label{tab:ec}
  \end{center}
\end{figure}

\section{Elliptic normal bases}\label{section:fast}

In this section, we assume that we are in the situation
of Section~\ref{section:ellipticbasis}.
So $E$ is an elliptic curve over a finite field
$\bK=\Fq$ and  $d\geqslant 2$ is  an integer.
Let $t\in E(\Fq)[d]$ be a rational point of order
$d$. We call $T$ the group generated by $t$. 
Let $\phi: E \rightarrow E$
be the Frobenius endomorphism. 
Let $b\in E(\bar \bK)$  be a point such that $\phi(b)=b+t$.
So, $b$   belongs to $E(\bL)$ where $\bL$ is the degree
$d$ extension of $\bK$. 
We denote by $E'$ the quotient $E/T$ and by $I: E\rightarrow E'$
the quotient isogeny. We also assume $db\not =O\in E$.
We set $a=I(b)$ and check $a\in E'(\Fq)$.
We further assume there  exists
one point $R$ in $E(\Fq)$ such that  $dR\not =0$.
\medskip

 We construct a normal basis for $\bL$, the degree $d=\# T$ extension 
of $\bK$. In this basis, the product  of two
elements can be computed at the expense of $5$ convolution products between
vectors of dimension $d$. Such  bases may  be preferred to the ones
constructed in Section~\ref{section:ellipticbasis} when $d$ is large enough,
depending on the implementation context.

\subsection{The elliptic normal basis $\Theta$}

We start with a lemma
concerning the sum
$\sum_{k\in \ZZ/d\ZZ} u_{kt,(k+1)t}$.
\begin{lemma}
The sum $\sum_{k\in \ZZ/d\ZZ} u_{kt,(k+1)t} $ 
is a constant $\cgot \in \bK$. If the characteristic $p$ 
of $\bK$ divides the degree $d$, then $\cgot \neq 0$.
\end{lemma}

\begin{myproof}
The sum $\sum_{k\in \ZZ/d\ZZ} u_{kt,(k+1)t} $ is invariant
by translations in $T$. So it can be seen
as a function on $E'=E/T$. As such, it has no more than one
pole. Therefore it is constant.

Assume now $p$ divides $d$ and $\sum_{k\in \ZZ/d\ZZ} u_{kt,(k+1)t} = 0$. 
The sum $\sum_{k\in \ZZ/d\ZZ} ku_{kt,(k+1)t}$ is thus 
invariant
by translations in $T$. So it can be seen
as a function on $E'=E/T$. As such, it has no more than one
pole. Therefore it is constant. However, seen as a function
on $E$,  this sum $\sum_{k\in \ZZ/d\ZZ} ku_{kt,(k+1)t}$  has a pole at $O$.
A contradiction.
\end{myproof}\\
So  at least one  of the two following conditions holds: either $d$ is
prime to $p$ or $\cgot \not =0$.
In any case,  there exist two scalars $\agot \not = 0$ and 
$\bgot$ in $\bK$
such that $\agot\cgot+d\bgot=1$.
For $k\in \ZZ/d\ZZ$ we set
$u_k=\agot u_{kt,(k+1)t}+\bgot$ and $x_k=x_{kt}$.
\medskip

We denote by $\Theta$
the system $(\theta_k)_{k \in \ZZ/d\ZZ}$
defined as $\theta_k=u_k(b)$. 
We have
$\sum_{k\in \ZZ/d\ZZ} \theta_k=1\in \bK$.
and  $\Phi(\theta_k)=\theta_{k-1}$.
\begin{lemma}
With the above notation, the system 
$(u_0, u_1, \ldots , u_{d-1})$ is
 a  basis of  
$$\cL=\cL(\sum_{k\in \ZZ/d\ZZ} [kt]).$$
The system 
$\Theta = (\theta_0, \theta_1, \ldots, \theta_{d-1})$ is
 a $\bK$ basis of  $\bL$.
\end{lemma}
\begin{myproof}
Indeed, let the  $\lambda_k$ for $k\in \ZZ/d\ZZ$
be scalars in $\bK$ such that $\sum_{k\in \ZZ/d\ZZ}\lambda_k\theta_k=0$.
The function $f=\sum_{k\in \ZZ/d\ZZ}\lambda_ku_k$ cancels
at $b$ and also at all its $d$ conjugates over $\bK$ (because $f$ is
defined over $\bK$). But $f$ has no more than $d$ poles (the points
in $T$). If $f$ is non-zero, its  divisor is
 $(f)_0-(f)_\infty$ with $(f)_0=\sum_{t\in T} [b+t]$
and $(f)_\infty=\sum_{t\in T}[t]$. We deduce 
$d\times b$ is zero in $E$. But this is impossible by hypothesis.
So $f$ is constant equal to zero.
This implies all $\lambda_k$'s are equal (look at poles). Since the sum of
all $\theta_k$'s is non-zero, this implies that
 all $\lambda_k$'s are null. 
\end{myproof}\\
We call such a basis as $\Theta$ an {\it elliptic normal basis}.
\medskip

If $k, l \in \ZZ/d\ZZ$ and $k\not = l,l+1,l-1 \bmod d$, then
\begin{displaymath}
u_ku_l\in \cL
\end{displaymath} 
where $\cL=\cL(\sum_{k\in \ZZ/d\ZZ} [kt])$ is the $\bK$-vector space
generated 
by all  $u_m$ for $m\in \ZZ/d\ZZ$.
Further
\begin{displaymath}
u_{k-1}u_{k} +\agot^2 x_k \in \cL
\text{ and }
u_{k}^2-\agot^2x_k-\agot^2x_{k+1}\in \cL\,.
\end{displaymath} 
So if $(\alpha_k)_{0\leqslant k  \leqslant d-1}$
and $(\beta_k)_{0\leqslant k \leqslant d-1}$ are two vectors in $\bK^d$, we have
\begin{eqnarray}\label{eq:mulprinc}
(\sum_{k} \alpha_ku_{k})(\sum_{k}
\beta_ku_{k})&=&\agot^2
\sum_k\alpha _k\beta_k(x_k+x_{k+1})-
\agot^2\sum_k \alpha_{k-1}\beta_{k}x_{k} - 
\agot^2\sum_k
\beta_{k-1}\alpha_{k}x_{k}\bmod \cL \nonumber \\
&=&\agot^2 \sum_k(\alpha_{k}-\alpha_{k-1})(\beta_k-\beta_{k-1})x_k\bmod \cL.
\end{eqnarray}

\paragraph{Example.}
Let us continue the example of section~\ref{section:ellipticbasis},
\textit{i.e.}  $\bK=\FF_{7}$ and $d=5$.  We find
\begin{displaymath}
  (u_{kt,(k+1)t})_k =
  \left(
{\frac {5\,y+3}{x+4}},
{\frac {5\,y+3\,{x}^{2}+4}{{x}^{2}+5}},
 \frac{4}{x+3},
{\frac {y\,(2\,x+8)+3\,{x}^{3}+15\,x}{\left( {x}^{2}+5 \right)  \left( x+4
\right) }},
{\frac {2\,y+2\,x+6}{x+4}}
\right),
\end{displaymath}
so that $\cgot=3$, $\agot = 5$, $\bgot = 0$, and
\begin{displaymath}
  \Theta=(\tau^{8083}, \tau^{13159}, \tau^{16285}, \tau^{9529},
\tau^{6163})\,.
\end{displaymath}

\subsection{Change of coordinates}\label{para:changecoordinates}

Thanks to Equation~(\ref{eq:somme}), the $\theta$'s  can be given in the basis $(\omega_k)_k$ as
\begin{displaymath}
  \theta_k = \left\{
    \begin{array}{ll}
    \agot \omega_{1} +\bgot\omega_0 & \text{if }k = 0,\\
    -\agot \omega_{-1}-a_1\agot \omega_0  +\bgot \omega_0  & \text{if }k = d-1,\\
    \agot \omega_{k+1}- \agot \omega_{k}+\agot\Gamma_{k,k+1}\,\omega_0   +\bgot\omega_0& \text{otherwise.}
  \end{array}\right.
\end{displaymath}
Inversely, we set  $\lambda_k=\sum_{i=1}^{k} \Gamma_{i, i+1}$ and we  
observe that $\cgot = \lambda_{d-2} -a_1$. We obtain
\begin{displaymath}
  \omega_k = \left\{
    \begin{array}{ll}
    \displaystyle \sum_{i=0}^{d-1}\theta_i  &
    \text{if }k = 0,\\
    \agot^{-1}\theta_0  - \bgot \agot^{-1}\sum_{i=0}^{d-1}\theta_i& \text{if }k = 1,\\
    \displaystyle -\agot^{-1}\theta_{-1}+(\bgot\agot^{-1}-a_1) \sum_{i=0}^{d-1}\theta_i& \text{if }k = -1,\\
    \displaystyle
\agot^{-1}\sum_{i=0}^{k-1}\theta_i-(k\bgot\agot^{-1}  +\lambda_{k-1}) \sum_{i=0}^{d-1}\theta_i & \text{otherwise.}
  \end{array}\right.
\end{displaymath}
This shows that one can compute the change of variable from $\Omega$
to $\Theta$, and back, at the expense of $O(d)$ operations in $\bK$.

\subsection{Complexities}

We exhibit an algorithm with quasi-linear complexity to multiply two
elements given in an elliptic normal basis. As often with FFT-like
algorithms, it consists in evaluations and
interpolations.\medskip

\noindent
\textbf{Notation.}
If  $\va = (\alpha_i)_{0\leqslant i\leqslant d -1}$ and 
$\vb = (\beta_i)_{0\leqslant i\leqslant d-1}$ are two  vectors of length $d$
we denote by  $\va \star_j \vb = \sum_i \alpha_i\beta_{j-i}$
the $j$-th component of the convolution product.
We denote by $\sigma(\va)=(\alpha_{i-1})_{i}$ the cyclic
shift of $\va$.
We denote by $\va \diamond \vb
=(\alpha_i\beta_i)_i$ the component-wise
product and by
$\va \star \vb=(\va \star_i \vb)_i$ the convolution product.

\subsubsection{Reduction  }

Given a linear combination of the $\xi$'s we may want to
reduce it: express it as a linear combination of
the $\theta$'s. \medskip

Let $\viota=(\iota_i)_{0\leqslant i\leqslant d-1}$ be the
  vector  in $\bK^d$ such that 
$\xi_0=\sum_{0\leqslant k\leqslant d-1} \iota_k \theta_k$.
$$\xi_i=\Phi^{-i}(\xi_0)=
\sum_{0\leqslant k \leqslant d-1}\iota_k\Phi^{-i}(\theta_k)=
\sum_{0\leqslant k \leqslant d-1}\iota_k\theta_{k+i}=
\sum_{0\leqslant k \leqslant d-1}\iota_{k-i}\theta_{k}.
$$
Let $\va = (\alpha_i)_{0\leqslant i\leqslant d -1}$ and 
$\vb = (\beta_j)_{0\leqslant j\leqslant d-1}$ be vectors in $\bK^d$ such 
that 
$$\sum_{0\leqslant i\leqslant d-1}\alpha_i\xi_i=\sum_{0\leqslant j\leqslant d-1}\beta_j\theta_j.$$
We want to express the $\beta_j$'s as linear expressions in
the $\alpha_i$'s.
\begin{eqnarray}\label{eq:fastreduc}
\sum_{0\leqslant i\leqslant d-1}\alpha_i\xi_i&=&\nonumber
\sum_{0\leqslant i\leqslant d-1}\alpha_i \sum_{0\leqslant k \leqslant d-1}\iota_{k-i}\theta_{k}\\
&=&\sum_k\theta_k\sum_{i}\alpha_i\iota_{k-i}=
\sum_k(\viota\star_k\va )\theta_k.
\end{eqnarray}
We deduce
$\vb=\viota \star \va.$
So $\vb$ is  the convolution  product of 
$\viota $ and $\va$.

\subsubsection{Evaluation}

Let $(\alpha_i)_{0\leqslant i\leqslant d-1}$ be scalars in $\bK$. 
Let $R\in E(\bK)-E[d]$ be a $\bK$-rational point
on $E$ such that $dR\not =0$.
\medskip

We want to evaluate $f=\sum_{0\leqslant i \leqslant d-1} \alpha_ix_i$
at all $R+jt$ for $0\leqslant j\leqslant d-1$. We set
$\beta_j=f(R+jt)$.
We have
$$\beta_j=\sum_{0\leqslant i \leqslant d-1} \alpha_ix_i(R+jt)
=\sum_{0\leqslant i \leqslant d-1} \alpha_ix_0(R+(j-i)t)=
\va \star_j \vxr $$
where $\vxr=(x_0(R+kt))_{0\leqslant k \leqslant d-1}$.
So,
\begin{equation*}
\vb = \vxr\star \va.
\end{equation*}

Similarly, we want
to evaluate $f=\sum_{0\leqslant i \leqslant d-1} \alpha_iu_i$
at all $R+jt$ for $0\leqslant j\leqslant d-1$. We set
$\beta_j=f(R+jt)$. 
We have
$$\beta_j=\sum_{0\leqslant i \leqslant d-1} \alpha_iu_i(R+jt)
=\sum_{0\leqslant i \leqslant d-1} \alpha_iu_0(R+(j-i)t)=
\va \star_j \vur $$
where $\vur=(u_0(R+kt))_{0\leqslant k \leqslant d-1}$.
So,
\begin{equation}\label{eq:fastevalII}
\vb = \vur\star \va.
\end{equation}

\subsubsection{Interpolation}

Let $R\in E(\bK)-E[d]$ be a $\bK$-rational point
on $E$ such that  $dR\not = 0$. 
The evaluation map $f\mapsto (f(R+jt))_{0\leqslant j\leqslant d-1}$
is a bijection from $\cL$ onto $\bK^d$.
\medskip

Given the $\beta_j=f(R+jt)$ we want to compute the
$\alpha_i$ such that  $f=\sum_{0\leqslant i \leqslant d-1} \alpha_iu_i$.
Since $\vb = \vur\star \va$ we just need to compute once
for all the inverse  $\vcur$ of $\vur$
for the convolution product. This
inverse exists because the evaluation map is bijective.

\subsubsection{Multiplication}\label{para:multipfast}

Let $\va=(\alpha_i)_{0\leqslant i\leqslant d-1}$ and 
$\vb = (\beta_i)_{0\leqslant i\leqslant d-1}$ be two vectors in $\bK^d$.
We want to multiply $\sum_{i}\alpha_i\theta_i$
and $\sum_i \beta_i\theta_i$.\medskip

We define four functions on $E$,
\begin{eqnarray*}
A&=&\sum_i \alpha_i u_i\,,\ B = \sum_i \beta_iu_i\,,\\
C&=& \agot^2 \sum_i (\alpha_i-\alpha_{i-1})(\beta_i-\beta_{i-1})x_i\,,\\
D&=&AB-C\,.
\end{eqnarray*}
The product 
we want to  compute is $A(b)B(b)=
C(b)+D(b)$.

From Equation~(\ref{eq:mulprinc}), we deduce
that $D$ is in $\cL$.
From  Equation~(\ref{eq:fastreduc}), we deduce that the
coordinates in $\Theta$ of $C(b)$
are given by the vector
$$\viota\star \left( \agot^2 (\va -\sigma(\va))\diamond (\vb-\sigma(\vb))\right).$$
According to
Equation~(\ref{eq:fastevalII}),
the evaluation of $A$ at the points $(R+jt)_j$ is given
by the vector
$\vur \star \va$. The evaluation
at these points of $D$ is
$(\vur\star \va)\diamond (\vur\star \vb)-\vxr\star 
(\agot^2(\va -\sigma(\va))\diamond (\vb-\sigma(\vb)))$.
If we $\star$ multiply this late  vector on the left
by $\vcur$ we obtain the coordinates of $D$
in the basis $(u_0,\ldots,u_{d-1})$. These are also the coordinates
of $D(b)$ in the basis $\Theta$.

Altogether, we have proved what follows.
\begin{lemma}
The multiplication tensor for  normal
elliptic  bases of type $\Theta$ is
\begin{multline*}
(\agot^2\viota)\star \left( (\va -\sigma(\va))\diamond (\vb-\sigma(\vb))
\right)+\\
\vcur \star \left((\vur\star \va)\diamond (\vur\star \vb)-(\agot^2\vxr)\star 
\left((\va -\sigma(\va))\diamond (\vb-\sigma(\vb))\right)\right)
\end{multline*}

It consists in $5$ convolution products, $2$ component-wise products,
 $1$ addition and $3$  subtractions between vectors of size
$d$, the degree of the extension. 

\end{lemma}

Note that convolution products can
be computed at the expense of $O(d\log d \log |\log d|)$ operations 
in $\bK$ using algorithms due to Sch\"onhage and Strassen \cite{SchonStras},
Sch\"onhage \cite{Schon}, and Cantor and Kaltofen \cite{CK}.

Note also that it is standard to use
elliptic curves (and even curves of higher genera) to bound
the bilinear complexity of multiplication. One should mention
in particular work by Chudnowsky \cite{Chud},  Shokrollahi \cite{Shokro}, 
Ballet \cite{Ballet}, Chaumine \cite{Chaum}. The tensor we produce here is not competitive
with theirs
from the point of view of bilinear complexity. But this tensor
 is symmetric 
enough to allow fast application of the Frobenius automorphism.

\paragraph{Example.} In the setting of the examples of
Section~\ref{section:ellipticbasis} and Section~\ref{section:fast},
\textit{i.e.}  $\bK=\FF_{7}$ and $d=5$, we first precompute, with
$R=(1,2)$ a point of order 10 on $E$,
\begin{displaymath}
  \viota=(0, 5, 5, 1, 0),\ \vur = (4, 1, 5, 1, 4),\ \vcur=(2, 2, 0, 4,
  0)\text{ and }\vxr=(1, 5, 5, 1, 2)\,.
\end{displaymath}
Now, we are going to multiply $\sum_{i}\alpha_i\theta_i$ and $\sum_i
\beta_i\theta_i$ with
\begin{math}
  \va=(6, 3, 6, 1, 2)\text{ and }\vb=(2, 6, 6, 4, 2)\,.
\end{math}
We first easily find $\va -\sigma(\va)=(4, 4, 3, 2, 1)$,
$\vb-\sigma(\vb)=(0, 4, 0, 5, 5)$
and thus $(\va -\sigma(\va))\diamond (\vb-\sigma(\vb))=(0, 2, 0, 3,
5)$\,.

Therefore,
\begin{eqnarray*}
  (\agot^2 \viota)\star \left( (\va -\sigma(\va))\diamond
(\vb-\sigma(\vb)) \right) &=&
  (6, 0, 4, 5, 5)\,,\\
  (\vur\star \va)\diamond (\vur\star \vb) &=&(0, 4, 0, 3, 0)\,,\\
  (\agot^2 \vxr)\star \left((\va -\sigma(\va))\diamond
(\vb-\sigma(\vb))\right) &=&(1, 1, 0, 1, 4)\,.
\end{eqnarray*}
It remains to compute
\begin{displaymath}
  \vcur \star ((\vur\star \va)\diamond (\vur\star \vb)-(\agot^2 \vxr)\star
  ((\va -\sigma(\va))\diamond (\vb-\sigma(\vb)))) = (4, 5, 4, 0, 1)\,,
\end{displaymath}
and finally, we obtain
\begin{displaymath}
  (\sum_{i}\alpha_i\theta_i) \times (\sum_i \beta_i\theta_i) =
3\,\theta_0+ 5\,\theta_1+ 1\,\theta_2+ 5\,\theta_3+ 6\,\theta_4\,.
\end{displaymath}

\section{Beyond Gauss periods}\label{section:beyond}

Complexity estimates in Subsection~\ref{para:complex}
and Subsection~\ref{para:multipfast}
suggest that an elliptic basis may be preferred to standard
normal basis.
\medskip

In this section we first show that the main  condition
for the existence of an elliptic basis is that the degree should
not be too large. This is explained  in Subsection~\ref{para:borne}.
If this condition is not fulfilled, we may translate the field
extension along a small auxiliary base change. This is explained 
in Subsection~\ref{para:trans}.
We recall in Subsection~\ref{para:inv} that 
 fast inversion using Lagrange's theorem and
addition chains is possible
 in the context of elliptic normal bases.
In  Subsection~\ref{para:fastinv} we associate
a well chosen
polynomial  basis  to  any elliptic basis. We explain 
how to fast  change coordinates between either
bases.
This gives a quasi-linear division algorithm
for  elliptic
bases.

\subsection{Existence conditions for elliptic bases}\label{para:borne}

Let $q$ be a power of a prime $p$.
Given a finite field $\Fq$ and an integer $d\geqslant 2$,  we want
to construct an elliptic basis for the degree $d$ extension 
of $\Fq$.

We first need some easy properties of the $d_q$
(cf. Definition~\ref{definition:dq}).
\begin{lemma}\label{lemma:extdq} 
Let $p$ be a prime and $q$ a power of $p$. Let $d\geqslant 2$ be an integer.
\begin{itemize}
\item If $d$ is prime to $q-1$ then $d_q=d$.
\item If $q-1$ is squarefree then $d_q\leqslant d^3$.
\item In any case $d_q\leqslant d^2(q-1)^2$.
\item If $f\geqslant 1$ is  an integer  prime to $d\varphi(d)$ then
$d_{q^f}=d_q$.
\end{itemize}
\end{lemma}

We can now give a sufficient condition for the existence of an  
elliptic basis.  The necessary background about elliptic curves
over finite fields can be found in chapter 5 of Silverman's
book \cite{Silverman}.
\begin{lemma}\label{lemma:construc}
Let $p$ be a prime and $q$ a power of $p$. Let $d\ge 2$ be an integer.
We assume  that 
\begin{equation*}
d_q \leqslant 2\sqrt q\,.
\end{equation*}

Then, there exists an elliptic curve $E$ over $\Fq$, a point $t$ of order
$d$ in $E(\Fq)$ and a point $b$ in $E(\Fqb)$ such that
$\phi(b)=b+t$ and the order of $b$ is a multiple of
$d^2$. In particular $db\not =0$. 
\end{lemma}

\begin{myproof}

There are at least too consecutive multiples of $d_q$ in the interval
$[q+1-2\sqrt q, q+1+2\sqrt q]$. One of them is not congruent to 
$1$ modulo $p$. We call $M=\lambda d_q$ this integer
and we set  $\tgot =q+1-M$ and $\Delta=\tgot^2-4q$. 
Let $\cO$ be the maximal order in $\QQ(\sqrt{\Delta})$.  
There exists
an ordinary elliptic curve $E$ over $\Fq$ such that $E$ has $M$
 points over $\Fq$ and $\End(E)=\cO$\,.  Let $\ell$ be a prime
divisor of $d$. We set $e_\ell=v_\ell(d)$.
\bigskip 

\noindent
\underline{\it Assume first that $\ell$ is prime to $q-1$}.
\smallskip

\noindent
It cannot 
divide both $q+1-\tgot$
and $\tgot^2-4q$.
So $\ell$ is prime to $\tgot ^2-4q$ and is unramified in $\ZZ[\phi]$
and in $\End(E)$.
If $\ell$ were inert, it would divide both $\phi-1$ and its conjugate $\bar
\phi -1$ and also the trace $\Tr(\phi-1)=\tgot-2$. Since $\ell$ divides 
$q+1-\tgot$
this would imply that $\ell$ divides $q-1$, a contradiction. 
So  $\ell$ splits
 in $\ZZ[\phi]$. Let  $\lgot=(\ell,\phi-1)$ be the ideal  in
$\End(E)$ above
 $\ell$ and 
containing $\phi-1$. This prime ideal divides
$\phi-1$ exactly $e$ times, where $e\ge e_\ell$
 is the valuation of $M$ at $\ell$. 
Let $\lambda$  be the unique  root
of $(X+1)^2-\tgot (X+1)+q$ in $\ZZ_\ell$ that is congruent
to $0$ modulo $\ell$. The $\ell$-adic valuation  of
$\lambda$ is $e$. 
The kernel of $\lgot^{e+e_\ell}$ is cyclic of order
$\ell^{e+e_\ell}$. The Frobenius $\phi$ acts on this group
as multiplication by $1+\lambda$.
Let $b_\ell$ be a generator of this group. We set
$t_\ell=\phi(b_\ell)-b_\ell$ and
we check that $t_\ell$ has order $\ell^{e_\ell}$ and is $\Fq$-rational.
Indeed $t_\ell$ is left invariant
by $\phi$ because $e\geqslant e_\ell$. 
\medskip 

\noindent
\underline{\it Assume now $\ell$ divides $q-1$}.
\smallskip

\noindent
So $v_\ell(M)\ge v_\ell(d_q)>2v_\ell(q-1)$.
We check
$$\tgot^2-4q=(q-1)^2+ M^2-2M(q+1)=(q-1)^2+O(\ell^s)$$
where $s=v_\ell(M)>2v_\ell(q-1)$
if $\ell $ is odd, and 
$s=v_\ell(M)+2>2v_\ell(q-1)+2$ if $\ell =2$.

We deduce $\tgot^2-4q$ is a square in $\ZZ_\ell$ and $\ell$
splits in $\End(E)$.
Let $\lambda_1$ and $\lambda_2$ be the two roots
of $(X+1)^2-\tgot (X+1)+q$ in $\ZZ_\ell$. Since 
$\lambda_1\lambda_2=q+1-\tgot=M$, one of these two roots
has $\ell$-adic valuation $\geqslant e_\ell$.
Assume for example $v_\ell(\lambda_1)=e_1\geqslant e_\ell$.
The $\ell^{e_1+e_\ell}$-torsion group
$E[\ell^{e_1+e_\ell}]$ has a cyclic subgroup
$V_1$ of order $\ell^{e_1+e_\ell}$ where
$\phi$ acts as multiplication by $1+\lambda_1$.

Let $b_\ell$ be a point of order $\ell^{e_1+e_\ell}$
in $V_1$. We set $t_\ell=\phi(b_\ell)-b_\ell=\lambda_1b_\ell$.
This is a point of order $\ell^{e_\ell}$. It is left invariant
by $\phi$ because $e_1\geqslant e_\ell$. So again $t_\ell$ 
is in $E[\ell^{e_\ell}](\Fq)$.
\medskip

\noindent
\underline{\it We now  patch all these points together}.
\smallskip

\noindent
We set $t=\sum_\ell t_\ell$ and $b=\sum_\ell b_\ell$. We have
$\phi(b)-b=t$ and $t$ has order $d$.
The order of the point $b$  is a multiple
of  $\prod_\ell \ell^{2e_\ell}=d^2$. In particular 
$db\not =0$. 
\end{myproof}

\begin{lemma}\label{lemma:construc2}
Let $p$ be a prime and $q$ a power of $p$. Let $d\ge 2$ be an integer.
We assume  that 
\begin{equation*}
d_q \leqslant \sqrt q\,.
\end{equation*}

Then, there exists an elliptic curve $E$ over $\Fq$, a point $t$ of order
$d$ in $E(\Fq)$ and a point $b$ in $E(\Fqb)$ such that
$\phi(b)=b+t$ and the order of $b$
is a multiple of $d^2$. In particular
$db\not =0$. There is also a point $R$ in $E(\Fq)$
that such that $dR\not =0$.
\end{lemma}
\begin{myproof}
We apply lemma \ref{lemma:construc} above to  $p$,
$q$ and  $d'=2d\leqslant  2\sqrt q$. We obtain an elliptic
curve $E$, a  point $t'$ of order $d'=2d$ 
in $E(\Fq)$ and a point $b'$ such that $\phi(b')=b'+t'$. We set
$t=2t'$,  $b=2b'$ and $R=t$ and we are done.
\end{myproof}

\subsection{Base change}\label{para:trans}

Let $q$ be a prime
power and let $d$ be an integer.
If $d$ is too large we may not be
able to construct an elliptic  basis for the  degree $d$ extension 
of $\Fq$. 
We try to  embed $\Fq$ into some small degree auxiliary extension
$\bK=\FQ$ with $Q=q^f$  then  construct an elliptic 
basis for the degree $d$ extension $\bL$ of $\bK$.
\medskip 
We shall need  the following lemma.
\begin{lemma}[Iwaniec]\label{lemma:Iwaniec}
  There exists a constant $K_{\mbox{\rm Iw}}\geqslant 1$ such that the
  following is true.

Let  $k\geqslant 2$ be an integer and let $p_1$, $p_2$, \ldots, $p_k$ 
be distinct prime integers. 
Let $\mu_i$ and $\mu_s$ be two integers with $\mu_s-\mu_i
\geqslant K_{\mbox{\rm Iw}}k^2(\log k)^2$. Let  $I$ be the interval $[\mu_i,
\mu_s]$. There is an integer $n$ in $I$  that is
prime to  every $p_i$ for $i\in \{1, 2, \ldots, k\}$.
\end{lemma}

This lemma is proven by  Iwaniec in \cite{Iwaniec}.

\medskip

The number of prime divisors of $d$ is $O(\log d)$.
We look for some integer $f$ such that
\begin{itemize}
\item $f$ is prime to $d\varphi(d)$\,,
\item $d_{q^f}
=d_q\leqslant q^{\frac{f}{2}}.$
\end{itemize}

From  Lemma~\ref{lemma:Iwaniec}, we find some 
$f$ that is 
$$O(\log_q d_q+(\log d)^2(\log(\log d))^2)=O((\log d)^2(\log(\log d))^2).$$  
In this context, we call $\Phi_q: \bar \FF_{q}\rightarrow \bar \FF_q$
 the absolute Frobenius of $\Fq$ and $\Phi_Q=\Phi_q^f$
the Frobenius of $\bK$. Once given an elliptic
basis for $\bL/\bK$, we can compute efficiently the action of
$\Phi_Q$. 
Let $F$ be an integer such that 
$1\leqslant F\leqslant d-1$ and $fF=1\bmod d$.  The restriction of $\Phi_Q^F$ to $\FF_{q^d}$ is 
$\Phi_q: \FF_{q^d} \rightarrow \FF_{q^d}$. We thus can compute efficiently
the Frobenius action on $\FF_{q^d}$ using the elliptic basis for $\bL/\bK$.

Elements in $\FF_{q^d}$ being represented and treated as elements in $\bL$, we
have a slight loss of efficiency: the size is multiplied by $f$. An element in $\FF_{q^d}$ is represented by $d\log Q$ bits instead of
$d\log q$.

\subsection{Inversion using Lagrange's theorem}\label{para:inv}

We have constructed models for finite fields where addition, multiplication
and Frobenius action can be quickly computed. We should worry now about
 inversion. 
\medskip

The inverse of $\alpha \in \FF_{q^d}$ can be computed
as $\alpha^{q^d-2}$  because of Lagrange Theorem.
This exponentiation can be done at the expense of $O(\log q+\log d)$
multiplications in $\FF_{q^d}$ using an addition chain for $d-1$
and another addition chain for $q-2$. This is 
\cite[Theorem 2]{Itoh} of Itoh and Tsujii generalized in 
\cite[Corollary 30]{Nocker} by von zur Gathen and N\"ocker.
The computation also requires $O(\log d)$  exponentiations
by powers of $q$. 

\subsection{Moving to a polynomial  basis and quasi-linear   inversion}\label{para:fastinv}

Using Lagrange's theorem for inversion is one of the
possible motivations for using normal bases but
it brings an extra  $\log q$
factor in the complexity. This may harm if $\log q$ is bigger
than any polynomial in $\log d$. So it makes sense to look for
an inversion algorithm that uses less than e.g.
$K d(\log d)^2\log |\log d|$ operations in $\Fq$
where $K$ does not depend on $d$ nor on $q$.

In this subsection we show that to  any elliptic basis one can associate
a polynomial r basis such that changing coordinates between either
bases can be done in  quasi-linear time.
This gives another algorithm for fast multiplication in elliptic
bases. More importantly, this allows fast division in elliptic
bases.

Let $\bK=\Fq$, $d$, $\bL$,  $E$, $t$ and $b$
 be as in the beginning of
Section \ref{section:fast}. We further assume 
$2db\not = 0$. This is guaranteed if we use
Lemma \ref{lemma:construc2}  and if $d\geqslant 3$.
The unitary polynomial

$$\Pi(x)= (x-x(b))(x-x(b+t))\cdots (x-x(b+(d-1)t))\in \bK[x]$$
\noindent is then irreducible.

In order to simplify the presentation, we shall assume in the following
that $d$ is odd.
There exist
a degree $(d+1)/2$ unitary polynomial $Y_1\in \bK[x]$
and a degree $\leqslant  (d-3)/2$ polynomial $Y_0 \in \bK[x]$
such that the function $Y_1(x)-yY_0(x)$ cancels at
$b$, $b+t$,\ldots, $b+(d-1)t$. Besides $Y_1$ and $Y_0$
are coprime and $Y_1(x)-yY_0(x)$ also cancels at $-db$.
We precompute
these two  polynomials.

We denote by $\cR\subset \bK(E)$  the ring of functions 
having no pole outside $\{O,t,2t,\ldots, (d-1)t\}$.
The ideal $\bgot\subset \cR$
of the closed subset $\{b,b+t, b+2t, \ldots, b+(d-1)t\}$
is generated by $\Pi(x)$ and $Y_1(x)-yY_0(x)$. 

The system $(1, u_{O,t}, \ldots, u_{O,(d-1)t})$ is a 
$\bK$-basis of $\cL_1=\cL(O+t+2t+\dots+(d-1)t)$ and reduction
modulo $\bgot$ (evaluation at $b$) defines a bijection
$\epsilon_1 : \cL_1\rightarrow \bK(b)=\bL$. The system 
$(1, u_{O,t}(b), \ldots, u_{O,(d-1)t}(b))$ is the elliptic basis
$\Omega$.

The system $(1, x,  x^2,  \ldots, x^{{d-1}})$ is free and
generates a subspace $\cL_2$ of $\cL((2d-2)O)$. Reduction
modulo $\bgot$ (evaluation at $b$) defines a bijection
$\epsilon_2 : \cL_2\rightarrow \bK(b)=\bL$. The system
$\Psi=(1,x(b),x(b)^2,\ldots,x(b)^{d-1})$ is a $\bK$-basis
of $\bL$.  This is a polynomial basis.

In order to change coordinates from $\Omega$
to $\Psi$ and back\footnote{Recall that
changing coordinates from $\Omega$ to $\Theta$ and back is done
in linear time as explained in praragraph \ref{para:changecoordinates}.}, 
we now
explain how to quickly evaluate  the bijections
$\epsilon_2^{-1}\circ \epsilon_1$ and $\epsilon_1^{-1}\circ \epsilon_2$.

\medskip 
\noindent
\underline{\it From $\Omega$ to $\Psi$}.
\smallskip

Recall we have set $\nu_k=x(kt)$ for $k\in \ZZ/d\ZZ$.
Equation (\ref{eq:gamma}) shows that there exist
constants $s_k = a_1x(kt)+a_3+y(kt)$ in 
$\bK$ such that for $1\leqslant  k \leqslant  d-1 $ 

$$u_{O,kt}=\frac{y+s_k}{x-\nu_k}.$$

Any function $f$ in $\cL_1$ is a combination 

$$f=\alpha_0+\sum_{1\leqslant  k\leqslant  d-1}\alpha_k\frac{y+s_k}{x-\nu_k}$$
\noindent with  $\alpha_k\in \Fq$ for $0\leqslant k \leqslant d-1$.
We set 

$$D(x)=\prod_{1\leqslant  k\leqslant  (d-1)/2}(x-\nu_k).$$

We can rewrite $f$ as ${(U(x)+yV(x))}/{D(x)}$
where $U(x)$ and $V(x)$ are polynomials in $\bK[x]$
with degree $\leqslant  \frac{d-3}{2}$. 

The numerator $U(x)+yV(x)$
can be computed
at the expense of $O(d(\log d)^2\log|\log d|)$  operations in $\Fq$ using
a divide and conquer algorithm.

Now the function $f$  is congruent modulo $\bgot$
to ${(U(x)+M(x)V(x))}/{D(x)}$. There exists a polynomial
$W(x)\in \bK[x]$ with degree $\leqslant  d-1$ that is congruent to the
later fraction modulo $\Pi(x)$. We compute it at the expense
of $O(d(\log d)^2 \log |\log d|)$ 
operations in $\Fq$ using standard fast  modular multiplication and inversion
 algorithms.
This polynomial $W(x)$ is nothing but 
$\epsilon_2^{-1}(\epsilon_1(f))$.

\medskip 

\noindent
\underline{\it From  $\Psi$ to $\Omega$}.
\smallskip

Conversely, let $W(x)\in \cL_2$  be a polynomial in $\bK[x]$
with degree $\leqslant  d-1$.
We look for a function  $f=\alpha_0+ \sum_{1\leqslant  k\leqslant  d-1}
  \alpha_k({y+s_k})/({x-\nu_k})$ in $\cL_1$
 that is congruent
to $W(x)$ modulo $\bgot$.

For $k\not = 0$
in $\ZZ/d\ZZ$ we set 

$$D_k(x)=\prod_{1\leqslant  l\leqslant  (d-1)/2, \,\,\,  l\not \equiv \pm k
\bmod d}(x-\nu_l)=D(x)/(x-\nu_k).$$

We assume we have precomputed the $D_k(\nu_k)$ for $1\leqslant  k\leqslant  (d-1)/2$
using   fast multipoint evaluation of the derivative $D'(x)$
at the expense of 
$O(d(\log d)^2 \log |\log d|)$ operations in $\Fq$.

We first compute a degree $\leqslant d-1$ polynomial 
 $N(x)$  that is congruent to $W(x)D(x)Y_0(x)$ modulo $\Pi(x)$.
This is done at the expense
of $O(d(\log d)^2 \log |\log d|)$ 
operations in $\Fq$ using a standard  fast modular multiplication  
and reduction  algorithm.


We have 

$$N(x)\equiv  
D(x)Y_0(x)f\equiv \alpha_0D(x)Y_0(x)+\sum_{1\leqslant k\leqslant d-1} 
\alpha_kD_k(x)(Y_1(x)+s_kY_0(x))\bmod \bgot.$$

The leftmost and rightmost terms in the above congruence
are polynomials in $x$ with degree $\leqslant d-1$. Therefore they
are equal.
Since $D_k=D_{-k}$,  we obtain 
$$N(x)=\alpha_0D(x)Y_0(x)+\sum_{1\leqslant k\leqslant (d-1)/2}
(\alpha_k(Y_1(x)+s_kY_0(x))+\alpha_{-k}(Y_1(x)+s_{-k}Y_0(x)))D_k(x).$$

We set 

\begin{equation}\label{eq:A0}
A_0(x)=\sum_{1\leqslant k\leqslant (d-1)/2}(\alpha_ks_k+\alpha_{-k}s_{-k})D_k
\text{ \, \, and  \, \,}
A_1(x)=\sum_{1\leqslant k\leqslant (d-1)/2}(\alpha_k+\alpha_{-k})D_k
\end{equation}
\noindent 
 and we obtain

$$N(x)=\alpha_0D(x)Y_0(x)+A_0(x)Y_0(x)+A_1(x)Y_1(x).$$

We now reduce this identity modulo $Y_1(x)$.
Let $\hat N(x)\in \bK[x]$ be a polynomial with degree 
$\leqslant (d-1)/2$ that is congruent to $N(x)/Y_0(x)$ modulo $Y_1(x)$.
We have $A_0(x)=\hat N(x)-\alpha_0 D(x)$ where $\alpha_0$ is the only
constant in $\bK$ such that $\hat N(x)-\alpha_0 D(x)$ has
degree $\leqslant (d-3)/2$.
Once we know $\alpha_0$ and $A_0(x)$ we set 
$A_1(x)=(N(x)-\alpha_0D(x)Y_0(x)-A_0(x)Y_0(x))/Y_1(x)$.

From equations (\ref{eq:A0}) we deduce

\begin{eqnarray*}
\alpha_ks_k+\alpha_{-k}s_{-k}&=&A_0(\nu_k)/D_k(\nu_k),\\
\alpha_k+\alpha_{-k}&=&A_1(\nu_k)/D_k(\nu_k).
\end{eqnarray*}

These pairs of  equations allow us to compute
all the $\alpha_k$ from the
$A_0(\nu_k)$, $A_1(\nu_k)$, and $D_k(\nu_k)$
at the expense of  $O(d)$
operations in $\bK$.
The $A_0(\nu_k)$ and  $A_1(\nu_k)$  are
computed using  a fast multipoint evaluation algorithm
at the expense of 
$O(d(\log d)^2 \log |\log d|)$ operations in $\Fq$.


\end{document}